\documentclass[12pt,a4paper]{amsart}
\usepackage{amscd}
\usepackage{mathrsfs}
\usepackage{amsfonts}
\usepackage{graphicx}
\usepackage{amsmath,amscd,stmaryrd,amssymb,array, amsfonts,amsmath,amssymb}
\usepackage{enumitem}
\usepackage{amsthm}
\usepackage{bm}

\setenumerate[1]{itemsep=0pt,partopsep=0pt,parsep=\parskip,topsep=5pt}
\setitemize[1]{itemsep=0pt,partopsep=0pt,parsep=\parskip,topsep=0pt}
\setdescription{itemsep=0pt,partopsep=0pt,parsep=\parskip,topsep=5pt}

\numberwithin{figure}{section}
 \numberwithin{equation}{section}

\newtheorem{theorem}{Theorem}[section]
\newtheorem{proposition}[theorem]{Proposition}
\newtheorem{definition}[theorem]{Definition}
\newtheorem{corollary}[theorem]{Corollary}
\newtheorem{lemma}[theorem]{Lemma}
\newtheorem{remark}[theorem]{Remark}
\newtheorem{example}[theorem]{Example}

\def\be{\begin{equation}}
\def\ee{\end{equation}}
\def\bes{\begin{equation*}}
\def\ees{\end{equation*}}
\def\bsp{\begin{split}}
\def\esp{\end{split}}

\def\ba{\begin{array}}
\def\ea{\end{array}}
\def\benu{\begin{enumerate}}
\def\eenu{\end{enumerate}}
\def\bt{\begin{theorem}}
\def\et{\end{theorem}}
\def\bp{\begin{proposition}}
\def\ep{\end{proposition}}
\def\bl{\begin{lemma}}
\def\el{\end{lemma}}
\def\br{\begin{remark}}
\def\er{\end{remark}}
\def\bd{\begin{definition}}
\def\ed{\end{definition}}

\def\W{\Omega}

\def\.{\cdot}

\def\~{\tilde}
\def\8{\infty}

\def\Vs{\vskip8pt}\def\vs{\vskip4pt}

\def\({\left(}\def\){\right)}

\begin{document}

\begin{center}
{\bf\Large
Multiple nodal solutions having shared componentwise nodal numbers
for coupled Schr\"{o}dinger equations}
\end{center}

\vs\centerline{Haoyu Li}  
\begin{center}
{\footnotesize
{Center for Applied Mathematics,  Tianjin University\\
          Tianjin 300072, P. R. China\\
{\em E-mail}:  hyli1994@hotmail.com}}
\end{center}
\vs\centerline{Zhi-Qiang Wang}
\begin{center}
{\footnotesize
{
Department of Mathematics and Statistics, Utah State University\\
           Logan, Utah 84322, USA\\
{\em E-mail}:  zhi-qiang.wang@usu.edu}}
\end{center}
\date{04.23.2019}

\Vs

{\footnotesize
{\bf Abstract.}
We investigate the structure of nodal solutions for coupled nonlinear Schr\"{o}dinger equations in the repulsive coupling regime.
Among other results, for the following coupled system of $N$ equations, we prove the existence of infinitely many nodal solutions which share the same
componentwise-prescribed nodal numbers
\begin{equation}\label{ab}
    \left\{
   \begin{array}{lr}
     -{\Delta}u_{j}+\lambda u_{j}=\mu u^{3}_{j}+\sum_{i\neq j}\beta u_{j}u_{i}^{2} \,\,\,\,\,\,\,  in\ \W ,\\
     u_{j}\in H_{0,r}^{1}(\W), \,\,\,\,\,\,\,\,j=1,\dots,N,
   \end{array}
   \right.
\end{equation}
where $\W$ is a radial domain in $\mathbb R^n$ for $n\leq 3$, $\lambda>0$, $\mu>0$, and $\beta <0$.
More precisely, let $p$ be a prime factor of $N$ and write $N=pB$. Suppose $\beta\leq-\frac{\mu}{p-1}$. Then for any given non-negative integers $P_{1},P_{2},\dots,P_{B}$, (\ref{ab}) has infinitely many solutions $(u_{1},\dots,u_{N})$ such that each of these solutions satisfies the same property: for $b=1,...,B$, $u_{pb-p+i}$ changes sign precisely $P_b$ times for $i=1,...,p$.
The result reveals the complex nature of the solution structure in the repulsive coupling regime due to componentwise segregation of solutions.
Our method is to combine a heat flow approach as deformation with a minimax construction of
the symmetric mountain pass theorem using a $\mathbb Z_p$ group action index.
Our method is robust, also allowing to give the existence of one solution without assuming any symmetry of the coupling.

 \Vs
{\bf Keywords:}
Multiple nodal solution; Componentwise-prescribed node; Coupled Schr\"{o}dinger equations.
\Vs {\bf 2010 MSC:} 35J47, 35J50, 35J55, 35K45}

\section{Introduction}

\subsection{Main Result}
In this paper, we consider the following coupled nonlinear Schr\"{o}dinger system of $N$ equations:
\begin{equation}\label{e:A111}
    \left\{
   \begin{array}{lr}
     -{\Delta}u_{j}+\lambda_{j} u_{j}=\mu_{j} u^{3}_{j}+\sum_{i=1, i\neq j}^N\beta_{ij} u_{j}u_{i}^{2} \,\,\,\,\,\,\,  in\ \W ,\\
     u_{j}\in H_{0,r}^{1}(\W), \,\,\,\,\,\,\,\,j=1,\dots,N,
   \end{array}
   \right.
\end{equation}
where $\W\subset \mathbb R^n$ is a radially symmetric domain, bounded or unbounded, $n\leq 3$, and the constants satisfy $\lambda_j>0$, $\mu_j>0$ for $j=1,...,N$, and $\beta_{ij}=\beta_{ji}$ for $i\neq j$. $H_{0,r}^{1}(\W)$ denotes the subspace of $H^{1}_{0}(\W)$ of radially symmetric functions.

To demonstrate the spirit of our results, we state the result in a special case first, where all $\lambda_j$ are equal $\lambda>0$,
all $\mu_j$ are equal $\mu>0$,
and all $\beta_{ij}$ are equal $\beta$ for $i\neq j$, i.e.,
\begin{equation}\label{e:All}
    \left\{
   \begin{array}{lr}
     -{\Delta}u_{j}+\lambda u_{j}=\mu u^{3}_{j}+\sum_{i\neq j}\beta u_{j}u_{i}^{2} \,\,\,\,\,\,\,  in\ \W ,\\
     u_{j}\in H_{0,r}^{1}(\W), \,\,\,\,\,\,\,\,j=1,\dots,N.
   \end{array}
   \right.
\end{equation}

\begin{theorem}\label{t:1}
  Let $p$ be a prime factor of $N$ and write $N=pB$. Suppose $\beta\leq-\frac{\mu}{p-1}$, Then for any given non-negative integers $P_{1},P_{2},\dots,P_{B}$, (\ref{e:All}) has infinitely many solutions $(u_{1},\dots,u_{N})$
  such that for $b=1,...,B$,
  $u_{pb-p+i}$ changes sign precisely $P_b$ times
  for $i=1,...,p$.
\end{theorem}

The result gives new insight into the structure of nodal solutions for coupled Schr\"odinger equations.
For a prescribed componentwise-node, we find infinitely many solutions which share the same number of nodal domains,  revealing more complexity of nodal solutions compared with the classical scalar field equation $-\Delta u +u = |u|^{p-2}u$ for which a long standing folklore has been the uniqueness of nodal solutions with a prescribed node. We say the solutions given above have componentwise-prescribed nodes.

Our method works in more general form than that of (\ref{e:All}).
Denote by $\mathcal B = (\beta_{ij})_{N\times N}$ the coefficient matrix involved on the right hand side of Problem (\ref{e:All}), where we denote $\beta_{ii}=\mu_{i}$.
We do not need to require the same values for $\lambda_j$, $\mu_j$ and $\beta_{ij}$. We denote the transformation of exchanging the $i$-th row and the $j$-th row of a matrix by $R_{ij}$, the $i$-th column and the $j$-th column by $C_{ij}$.

\begin{theorem}\label{t:main}
Let $p$ be a prime factor of $N$ and write $N=pB$. Assume the following four conditions hold.
\begin{itemize}
  \item [$(A)$] $\lambda_{pb-p+1}=\lambda_{pb-p+2}=\dots=\lambda_{pb}>0$ for $b=1,\dots,B$.
   \item [$(B)$] For $i,j=1,\dots,N$ and $i\neq j$, $\beta_{ij}=\beta_{ji}\leq0$ and $\mu_{j}>0$.
  \item [$(C)$] For $b=1\dots,B$, $\mathcal B = (\beta_{ij})_{N\times N}$ is invariant under the action of $$\prod_{i=1}^{p-1} C_{pb-p+i,pb-p+i+1}\circ R_{pb-p+i,pb-p+i+1}.$$

  \item [$(D)$] For $b=1,\dots,B$ and $pb-p+1\leq j\leq pb $, it holds $$\mu_{j}+\sum_{pb-p+1\leq i\leq pb\,;\,i\neq j}\beta_{ij}\leq 0.$$
\end{itemize}
Then for any given non-negative integers $P_{1},\dots,P_{B}$, the Problem (\ref{e:A111}) possesses infinitely many solutions $(u_{1},\dots,u_{N})$
such that for $b=1,...,B$,
  $u_{pb-p+i}$ changes sign precisely $P_b$ times
  for $i=1,...,p$.
\end{theorem}

\br
In the special case of \eqref{e:All}, we see (A) and (C) are satisfied readily, while (B) and (D) are satisfied under $\beta\leq-\frac{\mu}{p-1}<0$.
Thus Theorem \ref{t:1} follows.
\er

Our approach to study multiplicity of nodal solutions having the same componentwise nodes is to combine an associated parabolic flow serving as a descending flow of the variational problem with a minimax construction in the spirit of the symmetric mountain pass theorem via an $\mathbb Z_p$ index theory in the presence of invariant sets of the flow.
While for multiplicity of nodal solutions having the same nodal numbers, we need in an essential way the symmetry of the coupling coefficients, our method will be set up in a more general framework and also allows us to treat
the general case Problem (\ref{e:A111}) without a such symmetry. In this general setting we prove the existence of one solution with a prescribed node, and this gives a different proof of a result from \cite{LW2} where such a solution was given by gluing on Nehari manifold. In the present paper, we employ the corresponding parabolic flow as a tool for deformation of the variational problem, which is essential for establishing multiplicity results.

\begin{theorem}\label{t:exist}
Assume $\lambda_{j},\,\mu_{j}>0$ for $j=1\dots,N$. Then for any non-negative integers $P_{1},\dots,P_{N}$, there exists $b>0$ such that if $\beta_{ij}\leq b$ for all $i\neq j$, Problem (\ref{e:A111}) has a solution $(u_{1}\dots,u_{N})$ with the $j$-th component $u_{j}$ changing sign precisely $P_{j}$ times for $j=1\dots,N$.
\end{theorem}
We note that while for the multiplicity results we need the condition of negative coupling, for the existence of one solution we can allow a wider range of coupling here.


To make the symmetry condition in Theorem \ref{t:main} clear, we give three examples for the coupling coefficient matrix $\mathcal B = (\beta_{ij})_{N\times N}$ of Problem (\ref{e:A111}). The matrices are cut into blocks for suitable symmetry.

\begin{example}
For the case $N=4$ and $p=B=2$, the assumptions $(B-D)$ are satisfied in the following form
$$\left(
  \begin{array}{cc|cc}
    \mu_{1} & \beta_{1} & \beta_{3} & \beta_{3} \\
    \beta_{1} & \mu_{1} & \beta_{3} & \beta_{3} \\ \hline
    \beta_{3} & \beta_{3} & \mu_{2} & \beta_{2} \\
    \beta_{3} & \beta_{3} & \beta_{2} & \mu_{2} \\
  \end{array}
\right)$$
with $\beta_{i}\leq-\mu_{i}<0$ for $i=1,2$ and $\beta_{3}\leq0$. Assume that $\lambda_{1}=\lambda_{2}>0$ and $\lambda_{3}=\lambda_{4}>0$. Then given any two nonnegative integers $P_1, P_2$, there exist infinitely many solutions with first two components $u_1, u_2$ each having exactly $P_1$ simple zeros, and with the last two components $u_3, u_4$ each having exactly $P_2$ simple zeros.

If we set $N=6$, $p=2$ and $B=3$, the assumptions $(B-D)$ are satisfied in
$$\left(
    \begin{array}{cc|cc|cc}
      \mu_{1} & \beta_{1} & \beta_{4} & \beta_{4} & \beta_{5} & \beta_{5} \\
      \beta_{1} & \mu_{1} & \beta_{4} & \beta_{4} & \beta_{5} & \beta_{5} \\ \hline
      \beta_{4} & \beta_{4} & \mu_{2} & \beta_{2} & \beta_{6} & \beta_{6} \\
      \beta_{4} & \beta_{4} & \beta_{2} & \mu_{2} & \beta_{6} & \beta_{6} \\ \hline
      \beta_{5} & \beta_{5} & \beta_{6} & \beta_{6} & \mu_{3} & \beta_{3} \\
      \beta_{5} & \beta_{5} & \beta_{6} & \beta_{6} & \beta_{3} & \mu_{3} \\
    \end{array}
  \right)
$$
with $\beta_{i}\leq-\mu_{i}<0$ for $i=1,2,3$ and $\beta_{4},\beta_{5},\beta_{6}\leq0$.

If we set $N=6$, $p=3$ and $B=2$, the assumptions $(B-D)$ are satisfied in
$$\left(
    \begin{array}{ccc|ccc}
      \mu_{1} & \beta_{1} & \beta_{1} & \beta_{3} & \beta_{3} & \beta_{3} \\
      \beta_{1} & \mu_{1} & \beta_{1} & \beta_{3} & \beta_{3} & \beta_{3} \\
      \beta_{1} & \beta_{1} & \mu_{1} & \beta_{3} & \beta_{3} & \beta_{3} \\  \hline
      \beta_{3} & \beta_{3} & \beta_{3} & \mu_{2} & \beta_{2} & \beta_{2} \\
      \beta_{3} & \beta_{3} & \beta_{3} & \beta_{2} & \mu_{2} & \beta_{2} \\
      \beta_{3} & \beta_{3} & \beta_{3} & \beta_{2} & \beta_{2} & \mu_{2} \\
    \end{array}
  \right)
$$
with $\beta_{i}\leq-\frac{\mu_{i}}{2}<0$ for $i=1,2$ and $\beta_{3}\leq0$.
\end{example}

\subsection{Historical Remarks and the Idea of the Present Paper}
The nonlinear coupled elliptic system (\ref{e:A111}) has its theoretical root in Bose-Einstein condensates. The solutions of Problem (\ref{e:A111}) give rise to standing waves solutions of the time-dependent nonlinear coupled Schr\"{o}dinger system
\begin{equation}
    \left\{
   \begin{array}{lr}
     -i\partial_{t}\Phi_{j}-{\Delta}\Phi_{j}=\mu_{j}|\Phi_{j}|^2 \Phi_{j}+\sum_{i\neq j}\beta_{ij}\Phi_{j}|\Phi_{i}|^{2} \,\,\,\,\,\,\,  in\ \W ,\\
     \Phi_{j}(t,x)\in\mathbb{C}, \,\,\,\,\,\,\,\,j=1,\dots,N,
   \end{array}
   \right.
\end{equation}
for $j=1,\dots,N$ and $t>0$. In physics models, the parameters $\mu_{j}$ and $\beta_{ij}$ are the intraspecies and interspecies scattering lengths respectively. When $\beta_{ij}>0$, it is called the attractive case, when $\beta_{ij}<0$, it is called the repulsive case. In this paper, we mainly consider the repulsive case while small attractive coupling is also considered. \cite{AA, MS} is referred for more physics background.

 In recent years, a large number of mathematical results on Problem (\ref{e:A111}) have appeared, e.g., in \cite{AC, BDW, BW, DWW, LWei1, LWei2, LLW, LW1, LW, NTTV, S, TV, TW, WW} for studying various aspects of the problem such as existence theory for the attractive case and for the repulsive case, the bifurcation analysis, the synchronization and segregation for different coupling parameter regimes, and convergence and regularity of large couplings in the repulsive case etc.
 We refer to these papers for more references therein. In the repulsive coupling case, solutions tend to be segregated component-wisely creating more complex patterns of solutions.
 The application of variational methods to the coupled Schr\"{o}dinger systems mainly involves minimizing methods and minimax methods.
 The symmetric mountain-pass theorem has been well adopted for a large number of elliptic problems that goes back to the celebrated \cite{AR} by Ambrosetti-Rabinowitz. For Problem (\ref{e:A111}), the first difficulty is that there exist infinitely many so-called semi-trivial solutions (solutions with some components being zeros) so the system is degenerated to a system of smaller number of equations. In \cite{LW1,LW}, Liu and Wang proved the existence of infinitely many non-trivial (all components are non-zero) solutions to Problem (\ref{e:A111}) via invariant sets of descending flow and Nehari manifold method respectively. In \cite{DWW} and \cite{TW}, the authors proved multiplicity results of positive solutions to the special case Problem (\ref{e:All}) which possesses the componentwise permutation symmetry. This can be considered as a typical result for the repulsive case which shows distinct difference between a scalar field equation and a coupled nonlinear elliptic system since
 for the classical scalar field equation $-\Delta u +u = |u|^{p-2}u$ the uniqueness of positive solutions is well known (\cite{GNN, K}) and a folklore has been the uniqueness of nodal solutions with a prescribed node (\cite{AWY, Ta}).
In \cite{LLW}, the authors obtained a multiplicity result of solutions to Problem (\ref{e:A111}) in general domains with prescribed number of positive components and prescribed number of sign-changing components.
Recently, for radially symmetric domains, the existence of a nodal solution with componentwise prescribed number of nodes is obtained by Liu and Wang in \cite{LW2} via gluing on Nehari's method, extending the work for scalar equations (\cite{BWillem,Struwe}). More precisely, it is proved in \cite{LW2} that for any given nonnegative integers $P_1, ...,P_N$ there is a nodal solution $(u_1, ..., u_N)$ to Problem (\ref{e:A111}) such that $u_i$ has exactly $P_i$ simple zeros, $i=1,...,N$.

In the present paper,
our main concern and interest is that for a componentwise prescribed node
whether there are {\it multiple such solutions} sharing the given nodal number,
in particular whether there are {\it infinitely many such solutions?}
This is the main goal of our studies. Our result gives a
construction of infinitely many solutions sharing a given componentwise-prescribed
 node (Theorem \ref{t:1} and \ref{t:main}).



To deal with the sign-changing property of multiple solutions we will employ the heat flow of the corresponding coupled heat equations to Problem (\ref{e:A111}).
An important part of the present paper lies in the studies of the associated heat flow, including the existence and regularity results, the global existence and blow-up results, the non-increasing property of the sign-changing numbers along flow lines, the boundedness of trajectories and dynamical property of some invariant sets of the flow. We refer \cite{Amann,H,Lu,Q} for general discussions on the parabolic problems.
There have been a lot of works in the literature in which elliptic problems are solved with the help of heat flow methods. In \cite{CMT}, Conti, Merizzi and Terracini proved the existence of radial solutions with prescribed number of nodal domains to a scalar field equation. Utilizing the semilinear parabolic flow and the topological degree, they proved the result which was only treated by Nehari method before (\cite{BWillem,Struwe}). In \cite{Ch}, Chang established a variational framework and applied it to minimal surface problems. Quittner proved the existence and multiplicity of solutions of several semilinear elliptic problems and other dynamical properties by using parabolic flow in \cite{Q1, Q2, Q3}. In \cite{AB}, Ackermann and Bartsch developed the idea of superstable manifold and refined the symmetric mountain-pass theorem for sign-changing solutions (c.f. \cite[Section 2]{BWW}), which produced
 multiplicity results, nodal properties and order comparison results.
  More works on using parabolic flow to treat elliptic problems can be found in the references in these papers. However, there are few results on the coupled Schr\"{o}dinger systems using the heat flow. We mention \cite{WW} in which for two equations a comparison between components of positive solutions was obtained.
We will further develop the ideas in these papers by using heat flow as a tool of descending flow of our variational problem for Problem (\ref{e:A111}).
In fact, with the growth of the nonlinearity, a finer analysis on the global existence of the parabolic flow is also required. Combining the Cazenave-Lions interpolation (\cite{CL} and \cite{Ch}) and some estimate in \cite{Ch,Q1}, we can address that the growth of the nonlinearity is admissible to the global existence in dimensions $n\leq 3$. A finer analysis of the invariant sets requires the $H^1$-bounds of global solutions which we will prove in Section 2.4. We use a variant of the method in \cite{Q2003}, and we refer \cite{FL,Q2003} for more references on this topic.
Another important part of our work involves using some natural permutation symmetry in the coupling patterns.
 We will make use of the symmetry of the problem, that is, the problem is invariant under a $\mathbb Z_p$ group action of a cyclic permutation $\sigma$. With the heat flow serving as a deformation we will construct minimax critical values in the spirit of the symmetric mountain pass theorem via a $\mathbb Z_p$ index.
We need to build up special symmetric subsets of large $\mathbb Z_p$ index contained in the invariant sets of the flow.
Inspired by the approach  for scalar equations in \cite{CMT} our method is a sharper and symmetric variant of \cite{CMT} for coupled systems. To accomplish this, a certain combination of the methods in \cite{CMT,LW,TW} are needed.
While the idea of Nehari manifold was used in \cite{LW}, we will use the more natural ingredient, the boundary of the stable manifold of the origin, which has the advantage in keeping the non-increasing property of the sign-changing number along flow lines.

\subsection{The Structure of This Paper}
Section 2 mainly deals with the regularity and dynamical properties of the heat
flow of the corresponding heat equations, constructing various invariant sets of the flow.
We prove the existence result of Theorem \ref{t:exist} for the general system \eqref{e:A111} in Section 3, and this also will set up the stage for
the proof of the main result Theorem \ref{t:main}
in Section 4.
In Section 4, we give the proof of the multiplicity result Theorem \ref{t:main} from a minimax argument by constructing symmetric sets of large $\mathbb Z_p$
index inside various flow invariant sets on the boundary of the domain of attraction of the origin.

\section{Dynamical properties of the associated heat equations}

The parabolic flow associated with the elliptic system will be used as a mean of descending flow for the variational problems. We start by collecting some relevant results on existence
and regularity of the heat equations. Then we will develop some further estimates and construct some invariant sets of the flow which will be used in our proof later.
Let us fix some notations first.

We always use capital letters to represent vector valued functions and the corresponding lower case letters with subscript for their components.
For example, $U=(u_{1},\dots,,u_{N})$ and $V=(v_{1},\dots,v_{N})$.
A solution $U=(u_{1},\dots,,u_{N})$ to Problem (\ref{e:A111}) is called a non-trivial solution if and only if
$u_{j}\not\equiv0$ for any $j=1,\dots,N$. It is semi-trivial if and only if $U\not\equiv\theta$, where $\theta$ is the zero vector.

The norm of Lebesgue space $L^{p}(\W)$ is denoted by $|\cdot|_{p}$ and the norm of $H^{1}_{0}(\W)$ by $\|\cdot\|$.
For the product of spaces, such as $(L^{p}(\W))^N$ ($(H^{1}_{0}(\W))^N$),
we still use $|\cdot|_{p}$($\|\cdot\|$) to denote its norm. With no confusions, we sometimes omit the domain $\W$, the boundary condition and the radial condition and only denote the corresponding spaces by $L^p$, $H^{1}$, $H^{2}$ and $H^{s}$ for $s\in (1,2)$.

\subsection{Existence and Regularity Results of the Parabolic Flow}

Instead of the gradient flow, we will combine our variational structure with the following nonlinear coupled parabolic system:
\begin{equation}\label{e:A14}
    \left\{
   \begin{array}{lr}
     \frac{\partial}{\partial t}u_{j}-{\Delta}u_{j}+\lambda_{j} u_{j}=\mu_{j} u^{3}_{j}+\sum_{i\neq j}\beta_{ij} u_{j}u_{i}^{2} \,\,\,\,\,\,\,  in\ \W ,\\
     u_{j}(t,x)\in H_{0,r}^{s}(\W), \,\,\,\,\,\,\,\,j=1,\dots,N,\\
     u_{j}(0,x)=u_{0,j}(x)\in H_{0,r}^{s}(\W), \,\,\,\,\,\,\,\,j=1,\dots,N.
   \end{array}
   \right.
\end{equation}
whose equilibria are solutions to Problem (\ref{e:A111}).
Here, we require the coefficients $\lambda_j$'s, $\mu_j$'s and $\beta_{ij}$'s satisfy the conditions in Theorem \ref{t:exist} and Theorem \ref{t:main}
when we prove two theorems respectively.

By the notation $\eta^{t}(U)$ we denote a solution to the parabolic system with $U=(u_{0,1}, ..., u_{0,N})$ as its initial data. Sometimes,
  for the sake of simplicity, we also write $U(t)$.

  A special case of Problem (\ref{e:A14}) is of the form:
\begin{equation}\label{e:A12}
    \left\{
   \begin{array}{lr}
     \frac{\partial}{\partial t}u_{j}-{\Delta}u_{j}+\lambda u_{j}=\mu u^{3}_{j}+\sum_{i\neq j}\beta u_{j}u_{i}^{2} \,\,\,\,\,\,\,  in\ \W ,\\
     u_{j}(t,x)\in H_{0,r}^{s}(\W), \,\,\,\,\,\,\,\,j=1,\dots,N,\\
     u_{j}(0,x)=u_{0,j}(x)\in H_{0,r}^{s}(\W), \,\,\,\,\,\,\,\,j=1,\dots,N.
   \end{array}
   \right.
\end{equation}
It is obvious that an equilibrium point of Problem (\ref{e:A12}) is a solution to Problem (\ref{e:All}). Both of the parameters $s$ in
Problems (\ref{e:A14}) and (\ref{e:A12}) will be taken to be in $[1,2]$ depending upon the situation. Readers can find general theory of parabolic problems in \cite{Amann,DM,H,Lu,Q}. We will state a slightly more general result on the existence and regularity for the parabolic system (\ref{e:A14}) than we need in this paper.
Noticing that the spectrum of $-\Delta+\lambda$ is contained in $[\lambda,\infty)$, we conclude that the operator $-\Delta+\lambda$
is sectorial and, as a consequence, the existence and regularity results can be given. The results are stated and proved in terms of interpolation spaces $X_\alpha$ for $\alpha \in [0,1]$ (e.g., \cite{DM}).
We refer \cite{Amann,DM,H,Lu,Q} once
again for more information on sectorial operators and related properties on interpolation spaces. Note that with the range being $L^2$, the domain of the operator is $D(-\Delta+\lambda)=\{u\in H^{2}|\gamma_{2} u=0\}:=X_1$ (c.f. \cite[(4.6), (4.7)]{DM} or \cite{See}),
where $\gamma_{2}$ is the trace operator in $L^2=X_{0}$. Using the relation between these interpolation spaces and the Bessel-potential spaces $X_{s/2}=H_{0}^{s}(\W)$ for $\alpha\in\left[1,\frac{3}{2}\right)\cup\left(\frac{3}{2},2\right]$,
we will state the following theorem with the $H^s_0$ setting (c.f. \cite[Theorem 4.20]{DM}, \cite{PrS} and \cite{Amann87}). The following result for Problem (\ref{e:A14}) is useful in the present paper.
\begin{theorem}\label{t:2}
Let $s\in[1,2]\backslash\left\{\frac{3}{2}\right\}$ be fixed. If the initial value $U:=(u_{1},\dots,u_{N})\in(H^{s})^{N}$, there is a unique solution $\eta^{t}(U)=(u_{1}(t),\dots,u_{N}(t))$ to Problem (\ref{e:A14}) defined on its maximum interval $[0,T(U))$, satisfying
\begin{itemize}
  \item [(\uppercase\expandafter{\romannumeral1})]it holds that
  \begin{align}
   \eta^{t}(U)\in & C^1((0,T(U)),(L^2)^N) \cap C([0,T(U)),(H^s)^N);\nonumber
  \end{align}
  \item [(\uppercase\expandafter{\romannumeral2})] for any $U\in (H^{s})^{N}$ and any $\delta\in[0,T(U))$, there are positive constants $r,K$ such that for any $t\in[0,\delta]$
  $$\|U-V\|_{(H^{s})^{N}}<r\,\,\,\,\Rightarrow\,\,\,\,\|\eta^{t}(U)-\eta^{t}(V)\|_{(H^s)^N}\leq K\|U-V\|_{(H^{s})^{N}};$$
  \item [(\uppercase\expandafter{\romannumeral3})] the trivial solution $\theta\in (H^s)^N$ is asymptotically stable in $(H^s)^N$.
\end{itemize}
\end{theorem}



\br
Part $\mbox{(\uppercase\expandafter{\romannumeral1})}$ of this theorem is due to \cite[Theorem 15.3, Theorem 16.2]{DM} and $\mbox{(\uppercase\expandafter{\romannumeral2})}$ is of \cite[Proposition 16.8]{DM}. The assertion $\mbox{(\uppercase\expandafter{\romannumeral3})}$ is due to \cite[Theorem 5.1.1]{H}.
\er

\br
In following, we mainly use the result for $s=1$ and $s=2$.
\er


Notice that the theorem also holds true if we restrict the spaces to the case of radial symmetric functions. A similar regularity theory can be found in \cite{Lu}.

\subsection{Global Existence of the Solutions Starting on the Boundary of the Stable Manifold}
The propositions are modified versions of some results in \cite{CMT} and in \cite{Q1}. In this section, we always assume
$U(t)=(u_{1}(t),\dots,u_{N}(t))$ is a solution to Problem (\ref{e:A14}).
The energy of Problem (\ref{e:A111}) is the functional $J(U):(H_{r}^{1})^{N}\to\mathbb{R}$ defined by
\begin{align}
J(U)&:=J(u_{1},\dots,u_{N})\nonumber\\
&=\frac{1}{2}\sum_{j=1}^{N}\int|\nabla u_{j}|^{2}+\lambda_{j}|u_{j}|^{2}-\frac{1}{4}\sum_{j=1}^{N}\left(\int\mu_{j} u_{j}^4+\sum_{i\neq j}\int\beta_{ij} u_{i}^{2}u_{j}^{2}\right),\nonumber
\end{align}
which is a $C^{2}$ functional and satisfies the (PS) condition.

\begin{proposition}\label{p:L2partial_t}
For a solution $U(t)=(u_{1}(t),\dots,u_{N}(t))$, we have
$$\frac{\partial}{\partial t}J(U(t))=-\sum_{j=1}^{N}\int|\partial_{t}u_{j}|^2.$$
\end{proposition}
\noindent{\bf Proof.}
Note that
$$\frac{\partial}{\partial t}J(U)=\sum_{j=1}^{N}\nabla_{u_{j}}J(U)\partial_{t} u_{j}.$$
By a direct computation we have
\begin{align}
\nabla_{u_{j}}J(U)\partial_{t} u_{j}&=\int\nabla u_{j}\cdot\nabla\partial_{t}u_{j} +\lambda_{j} u_{j}\partial_{t}u_{j}-\int\mu_{j}  u_{j}^{3}\partial_{t}u_{j}+\sum_{i\neq j}\beta_{ij} u_{j}\partial_{t}u_{j}u_{i}^{2}\nonumber\\
                     &=\int \partial_{t}u_{j}\Bigg(-\Delta u_{j}+\lambda_{j} u_{j}-\mu_{j} u_{j}^{3}-\sum_{i\neq j}\beta_{ij}  u_{j}u_{i}^{2}\Bigg)=-\int|\partial_{t}u_{j}|^{2}.\nonumber
\end{align}
Then the proposition follows.
\begin{flushright}
\qedsymbol
\end{flushright}



\begin{corollary}\label{c:A} Let
\begin{align}
\mathcal{A}=\big\{U\in(H^{1})^{N}|T(U)=\infty\,\,\,\mbox{and}\,\,\lim_{t\to T(U)}\eta^{t}(U)=\theta\,\,\,\mbox{in}\,\,(H^1)^N\big\}.\nonumber
\end{align}
Then $\mathcal{A}$ is invariant under the heat flow and is open in $(H^{1})^{N}$.
\end{corollary}

This is a direct consequence of Theorem \ref{t:2}.


\begin{lemma}
$\partial\mathcal{A}$ is invariant under the heat flow and $\inf_{U\in\partial\mathcal{A}}J(U)\geq0$.
\end{lemma}
\noindent{\bf Proof.}
The continuity of the energy $J(U)$ with respect to $(H^1)^N$ norm implies the second part of this lemma. Now we prove the first part.

Due to the definition of $\mathcal{A}$, if $U\in\partial\mathcal{A}\subset (H^1)^N\backslash\mathcal{A}$, an immediate consequence is that $\eta^t(U)\subset(H^1)^N\backslash\mathcal{A}$. And if there is a $t_{0}\in(0,T(U))$ such that $\eta^{t_{0}}(U)\in(H^1)^N\backslash\overline{\mathcal{A}}$, due to Theorem \ref{t:2} and the openness of $\mathcal{A}$ in $(H^1)^N$, we can find a $V\in\mathcal{A}$ such that $\eta^{t_{0}}(V)\in(H^1)^N\backslash\overline{\mathcal{A}}$. Therefore, we address a contradiction. The above deduction implies that $\eta^{t}(U)\subset\overline{\mathcal{A}}\backslash\mathcal{A}=\partial\mathcal{A}$.

\begin{flushright}
\qedsymbol
\end{flushright}

Now, we prove that the flow with its initial data on the boundary of the stable manifold $\partial\mathcal{A}$ has $[0,\infty)$ as its maximal existence interval. Before that, let us prove a lemma under a more general condition.

\begin{lemma}
Suppose $\lim_{t\to T(U)}J(U)-J(\eta^{t}(U))\leq C<\infty$. Then $U(t)$ exists globally in $(H^1)^N$.
\end{lemma}
\noindent{\bf Proof.}
The proof makes use of some arguments from \cite{CL}, \cite[Lemma 1]{Ch} and \cite[Section 3]{Q1}. Since $\frac{\partial}{\partial t}J(\eta^t (U))=-\sum_{j=1}^{N}\int|\partial_{t}u_{j}|^{2}$, the condition in the lemma implies that
\begin{align}\label{ineq:PA1}
\sum_{j=1}^{N}\int_{0}^{t}\int|\partial_{t}u_{j}(s)|^{2}dxds=\big|J(\eta^{t}(U))-J(U)\big|\leq C.
\end{align}
Now we give the $L^2$-estimate. First we have
\begin{align}
\sum_{j=1}^{N}\int|u_{j}(t)|^2 dx&=\int_{0}^{t}\frac{d}{dt}\int|u_{j}(s)|^2 dxds+\sum_{j=1}^{N}\int|u_{j}(0)|^2 dx\nonumber\\
&=2\sum_{j=1}^{N}\int_{0}^{t}\int u_{j}\cdot\partial_{t}u_{j}dxds+\sum_{j=1}^{N}\int|u_{j}(0)|^2 dx\nonumber\\
&\leq C\Bigg(1+\sum_{j=1}^{N}\int_{0}^{t}\int|u_{j}(s)|^2 dxds\Bigg).\nonumber
\end{align}
Using the Gronwall's inequality, we have
\begin{align}\label{ineq:PA2}
\sum_{j=1}^{N}\int|u_{j}(t)|^2 dx\leq Ce^{Ct}.
\end{align}
Notice that
\begin{align}
&\;\;\;\;\sum_{j=1}^{N}\int|u_{j}(t)|^2 dx-\sum_{j=1}^{N}\int|u_{j}(0)|^2 dx\nonumber\\
&=2\sum_{j=1}^{N}\int_{0}^{t}\int u_{j}(s)\cdot\partial_{t}u_{j}(s)ds\nonumber\\
&=-8\int_{0}^{t}J(\eta^{s}(U))ds+2\sum_{j=1}^{N}\int_{0}^{t}\int|\nabla u_{j}(s)|^2 dxds.\nonumber
\end{align}
Therefore,
\begin{align}\label{ineq:PA3}
\int_{0}^{t}\int|\nabla u_{j}(s)|^2 dxds\leq Ce^{Ct}
\end{align}
follows immediately. Multiplying $u_{j}$ on both sides of the $j$-th equation of Problem (\ref{e:A14}), integrating over $\W$ and summing up with $j$, we obtain
$$\sum_{j=1}^{N}\int u_{j}\cdot\partial_{t}u_{j}+\sum_{j=1}^{N}\|u_{j}\|^{2}= \sum_{j=1}^{N}\int \mu_{j} u_{j}^{4}+\sum_{i\neq j}\beta_{j} u_{i}^{2}u_{j}^{2}.$$
Now we apply some methods in \cite{CL} and \cite{Ch}. For any $T>0$, we consider the norms on the time interval $[0,T]$. Due to the definition of the energy $J$ and (\ref{ineq:PA2}), we have
\begin{align}\label{inequality:important!!}
\sum_{j=1}^{N}\|u_{j}\|^2 &\leq 4J(U)+\sum_{j=1}^N \int u_{j}\cdot\partial_{t}u_{j}\nonumber\\
&\leq C+\Bigg(\sum_{j=1}^N |u_{j}|_{2}^{2}\Bigg)^{\frac{1}{2}}\cdot\Bigg(\sum_{j=1}^N |\partial_{t}u_{j}|_{2}^{2}\Bigg)^{\frac{1}{2}}\\
&\leq C+C\Bigg(\sum_{j=1}^N |\partial_{t}u_{j}|_{2}^{2}\Bigg)^{\frac{1}{2}}\nonumber
\end{align}
for $t\in[0,T]$. This implies that
\begin{align}
\int_{0}^{T}\big(\int|u_{j}(t)|^{2^{*}}dx\big)^{\frac{4}{2^*}}dt &\leq C\int_{0}^{T}\Bigg(\sum_{j=1}^{N}\|u_{j}(t)\|^{2}\Bigg)^2 dt\nonumber\\
&\leq C(T)+\sum_{j=1}^{N}\int_{0}^{T}\int|\partial_{t}u_{j}(t)|_{2}^{2}dt\leq C(T).\nonumber
\end{align}
 That is, $u_{j}\in L^{4}((0,T),L^{2^{*}} (\W))$, with $2^{*}=6$ for dimension 3. At the same time, the embedding $H^{1}\hookrightarrow L^{p}$ holds for any $p\geq 2$ for dimension 2. Therefore, we also have $u_{j}\in L^{4}((0,T),L^6(\W))$ for dimension 2. Notice that (\ref{ineq:PA1}) implies that $\partial_{t}u_{j}\in L^{2}((0,T),L^{2}(\W))$.

 Next we claim $u\in L^{\infty}((0,T),L^{\frac{18}{5}}(\W))$. We prove it in the following paragraph. Using the idea of  Cazenave-Lions interpolation (c.f. \cite{CL}, \cite{Ch}), we set $v=|u|^{3}$. Let us extend $u$ to $\W\times(0,3T)$ with compact support with respect to the variable $t$. Using the Newton-Leibniz formula and H\"{o}lder's inequality, we can compute that
\begin{align}
|u|_{\frac{18}{5}}&=|v|_{\frac{6}{5}}^{\frac{5}{6}\cdot\frac{1}{3}}=\Big|\int_{t}^{3T}|v_{t}|_{\frac{6}{5}}\Big|^{\frac{1}{3}} \leq C\Big|\int_{t}^{\infty}|u^{2}\cdot u_{t}|ds\Big|^{\frac{1}{3}}\nonumber\\
&\leq C\Big(\int_{0}^{3T}|u^{2}|_{3}|u_{t}|_{2}ds\Big)^{\frac{1}{3}}\nonumber
\end{align}
with
$$\frac{5}{6}=\frac{1}{3}+\frac{1}{2}.$$
Hence,
\begin{align}\label{ineq:CL1}
|u|_{\frac{18}{5}}&\leq C\Big(\int_{0}^{3T}|u|^{2}_{6}|u_{t}|_{2}ds\Big)^{\frac{1}{3}}\leq C\Big(\int_{0}^{3T}|u|^{4}_{6}ds\Big)^{\frac{1}{6}} \Big(\int_{0}^{3T}|u_{t}|^{2}_{l}ds\Big)^{\frac{1}{6}}\nonumber
\end{align}
with
$$1=\frac{1}{2}+\frac{1}{2}.$$
Therefore, we will have
\begin{align}
\sup_{t\in(0,T)}|u|_{\frac{18}{5}}&\leq C\Big(\int_{0}^{3T}|u|^{4}_{6}ds\Big)^{\frac{1}{6}}\Big(\int_{0}^{3T}|u_{t}|^{2}_{2}ds\Big)^{\frac{1}{6}}\nonumber\\
& \leq C(T)\nonumber
\end{align}
by using H\"{o}lder inequality with respect to the variable $t$. This implies that $u\in L^{\infty}((0,T),L^{\frac{18}{5}}(\W))$.

Multiplying $u_{j}^{3}$ on the both sides of the $j$-th equation and integrating over $\W$,
\begin{align}
\frac{d}{dt}|u^{2}_{j}|^{2}_{2}+\|u_{j}^{2}\|^{2} &\leq C\Big(\mu_{j}|u_{j}^{2}|_{3}^{3}+\sum_{i\neq j}\beta_{ij} u_{j}^{4}u_{i}^{2}\Big)\nonumber\\
&\leq C|u_{j}^{2}|_{3}^{3}.\nonumber
\end{align}
Using the interpolation inequality, we have
$$\frac{1}{3}=\frac{\frac{4}{7}}{6}+\frac{\frac{3}{7}}{\frac{9}{5}}$$
and
\begin{align}
\frac{d}{dt}|u^{2}_{j}|^{2}_{2}+\|u_{j}^{2}\|^{2} &\leq C|u_{j}^{2}|_{\frac{9}{5}}^{\frac{9}{7}}|u_{j}^{2}|_{6}^{\frac{12}{7}}\nonumber\\
&\leq C|u_{j}^{2}|_{\frac{9}{5}}^{\frac{9}{7}}\|u_{j}^{2}\|^{\frac{12}{7}}\nonumber\\
&\leq \frac{1}{2}\|u_{j}^{2}\|^{2}+C|u_{j}|_{\frac{18}{5}}^{18}\nonumber\\
&\leq \frac{1}{2}\|u_{j}^{2}\|^{2}+C.\nonumber
\end{align}
It follows that
$$|u_{j}^{2}(t)|^{2}_{2}\leq C\int_{0}^{T}dt+|u_{j}^{2}(0)|_{2}^{2}\leq C(T),$$
i.e. $u_{j}\in L^{\infty}((0,T), L^{4}(\W))$ for $j=1\dots,N$. Using the definition of the energy $J$,
$$\sum_{j=1}^{N}\|u_{j}\|^{2}\leq 2J(\eta^{t}(U))+\frac{1}{2}\sum_{j=1}^{N}\int\mu u_{j}^{4}+\sum_{i\neq j}\beta u_{i}^{2}u_{j}^{2} \leq C.$$
Hence $u_{j}\in L^{\infty}((0,T),H_{0}^{1}(\W))$.
Therefore, $T(U)=\infty$ since $T>0$ is arbitrary.

\begin{flushright}
\qedsymbol
\end{flushright}

\br\label{remark:globalexistence}
With the same conditions, we can address that $T(U)=\infty$, i.e. $\eta^t(U)\in C((0,\infty),(H^s)^N)$ for any $s\in (1,2]\setminus \{\frac{3}{2}\}$. To this end, we only need to use the formula of variation of constants and the fact that $\sum_{j=1}^{N}\Big(\mu_{j}u_{j}^{3}+\sum_{i\neq j}u_{i}^{2}u_{j}\Big)\in L^{\infty}((0,T),L^{2}(\W))$ (c.f. \cite[Proposition 7.1.8]{Lu}). Especially, with the same method, we can prove $\eta^t(U)\in C((0,\infty),(H^2)^N)$ when $U\in(H^2)^N$.
\er

The following corollary follows from the fact that $\inf_{\partial\mathcal{A}}J\geq 0$ and that the energy $J$ is non-increasing along the flow line.

\begin{corollary}
For any $U\in\partial\mathcal{A}$, $T(U)=\infty$.
\end{corollary}

\subsection{The $H^1$ bounds of the Solutions Starting on $\partial \mathcal{A}$}
We will borrow some ideas used by Quittner in \cite{Q2003}, which can be adapted for our situation (see also \cite{FL} for some related work).

\begin{lemma} \label{l2}
Let $U(t)$ be a global solution to (\ref{e:A14}) such that $\lim_{t\to\infty} J (U(t))= E_{1}$ is finite.
Then there is $C>0$ depending continuously upon the $L^2$-norm of the initial data, the initial energy $E_{0}:=J(U(0))$ and $E_{1}$, such that for any $t\geq0$, $|U(t)|_{2}\leq C$.
\end{lemma}
\noindent{\bf Proof.}
Suppose $0\leq t_{0}\leq t<+\infty$. Denote $\Phi(t)=\int_{t_{0}}^{t}|U(s)|_{2}^{2}ds$ and $E_{0}=J(U(0))$. Using the computation in previous subsection, we have
\begin{align}
\left |\sum_{j=1}^{N}|u_{j}(t)|_{2}^{2}-\sum_{j=1}^{N}|u_{j}(t_{0})|_{2}^{2} \right |&=2\left| \sum_{j=1}^{N}\int_{t_{0}}^{t}\int_{\W}u_{j}\partial_{t}u_{j}dxds \right|\nonumber\\
&\leq2\Big(\sum_{j=1}^{N}\int_{t_{0}}^{t}|\partial_{t}u_{j}|_{2}^2 ds\Big)^{\frac{1}{2}}\Big(\sum_{j=1}^{N}\int_{t_{0}}^{t}|u_{j}|_{2}^2 ds\Big)^{\frac{1}{2}}\nonumber\\
&\leq2\sqrt{E_{0}-E_{1}}\Phi(t)^{\frac{1}{2}},\nonumber
\end{align}
which gives
\begin{align}
\Phi'(t) \leq |U(t_{0})|_{2}^{2} + 2\sqrt{E_{0}-E_{1}} \Phi(t)^{\frac{1}{2}}.
\end{align}
Then we can compute that
\begin{align}
2\Big(\sqrt{\Phi(t)}-|U(t_{0})|_2\Big)'_{+}=\frac{\Phi'(t)}{\sqrt{\Phi(t)}}\chi_{\{\Phi>|U(t_{0})|_{2}\}}\leq|U(t_{0})|_{2}+2\sqrt{E_{0}-E_{1}}.\nonumber
\end{align}
This gives
\begin{align}
\sqrt{\Phi(t)}\leq|U(t_{0})|_2+\big(|U(t_{2})|_2+2\sqrt{E_{0}-E_{1}}\big)\frac{t-t_{0}}{2}.
\end{align}
Combining above deductions, we have
\begin{align}\label{inequ:L2bound3}
&\left |\sum_{j=1}^{N}|u_{j}(t)|_{2}^{2}-\sum_{j=1}^{N}|u_{j}(t_{0})|_{2}^{2}\right |\nonumber\\
\leq 2 &\sqrt{E_{0}-E_{1}}\Bigg(|U(t_{0})|_2 +\big(|U(t_{0})|_2+2\sqrt{E_{0}-E_{1}}\big)\frac{t-t_{0}}{2}\Bigg).
\end{align}
Set
\begin{align}
C_{1}=\frac{1}{\min_{j}\lambda_{j}}(9+8E_{0}+8E_{1})+2|U(0)|_{2}^{2}+81(E_{0}-E_{1})+3
\end{align}
We claim $|U(t)|_{2}^{2}\leq C_{1}$ for any $t\geq0$.
We prove the claim by contradiction. Suppose that there is a $\tau>0$ such that $|U(\tau)|_{2}^2 >C_{1}$.
First since $\sum_{j=1}^{N}\int_{0}^{\infty}|\partial_{t}u_{j}|_{2}^{2}ds=E_{0}-E_{1}<\infty$, we can find a sequence $t_{k}\to\infty$ such that $\nabla J(U(t_{k}))\to0$ and $J(U(t_k))\to E_{1}$ as $k\to\infty$. Thus $(U(t_k))\subset (H^1)^N$ is a (PS) sequence.
It is easy to check (e.g., \cite{W}) that $\|U_{k}\|^{2}\leq4(1+E_{1})$.
We also have
\begin{itemize}
  \item $|U(0)|_{2}^2 \leq |U(0)|_{2}^2 +1< \frac{C_{1}}{2}$;
  \item $|U(t_k)|_{2}^2 \leq\frac{4(1+E_{1})}{\min_{j}\lambda_{j}}<\frac{C_{1}}{2}$.
\end{itemize}
Let $k$ be the integer such that $\tau\in[t_{k-1},t_k]$ and, without loss of generality, let us assume $|U(\tau)|_{2}^{2}=\max_{[t_{k-1},t_k]}|U(s)|_{2}^{2}$. Then for any $t\in[\tau,\tau+1]$, applying (\ref{inequ:L2bound3}) and the fact that $C_{1}>81(E_{0}-E_{1})+1$,
\begin{align}
|U(t)|_{2}^{2}&\geq|U(\tau)|_{2}^{2}-2\sqrt{E_{0}-E_{1}}\Big(\frac{5}{2}|U(\tau)|_{2}+3\sqrt{E_{0}-E_{1}}\Big)\nonumber\\
&\geq|U(\tau)|_{2}^{2}-5\sqrt{E_{0}-E_{1}}|U(\tau)|_{2}-6(E_{0}-E_{1})\nonumber\\
&>\frac{|U(\tau)|_{2}^{2}}{2}>\frac{C_{1}}{2}>|U(t_k)|_{2}^{2}.\nonumber
\end{align}
This implies that $t_{k}\notin[\tau,\tau+1]$. Therefore $\tau+1<t_{k}$ and $\tau+1\in[t_{k-1},t_{k}]$. Consequently we have $|U(\tau+1)|_{2}\leq|U(\tau)|_{2}$. From above computation, we also have $|U(t)|_{2}^{2}\geq\frac{C_{1}}{2}$ for $t\in[\tau,\tau+1]$. And now, since $C_{1}>\frac{1+8(E_{0}-E_{1})}{\min_{j}\lambda_{j}}$, we have
\begin{align}
0&\geq|U(\tau+1)|_{2}^{2}-|U(\tau)|_{2}^{2}=2\sum_{j=1}^{N}\int_{\tau}^{\tau+1}\int_{\W}u_{j}\partial_{t}u_{j}dxds\nonumber\\
&=-8\int_{\tau}^{\tau+1}J(U(s))ds+2\int_{\tau}^{\tau+1}\sum_{j=1}^{N}\lambda_{j}|u_{j}|_{2}^{2}(s)ds\nonumber\\
&\geq-8E_{0}+2\min_{j}\lambda_{j}\int_{\tau}^{\tau+1}|U(s)|_{2}^{2}ds\nonumber\\
&\geq-8E_{0}+\min_{j}\lambda_{j}C_{1}>1,\nonumber
\end{align}
which is a contradiction. Hence we have $|U(t)|_{2}^{2}\leq C_{1}$ for any $t\geq0$.
\begin{flushright}
\qedsymbol
\end{flushright}

For a further discussion on the boundedness of trajectories of the flow in $(H^{1})^N$,
we need a result on the  maximal regularity of parabolic equation in \cite[Theorem \uppercase\expandafter{\romannumeral3}.4.10.7]{Amann}.
Spaces involving time will be used here, such as $L^{p}(I,X)$, $W^{1,p}(I,X)$, where $I$ is an interval and $X$ is a Banach space with a norm $||\cdot||_X$, and the norms are defined respectively by
      $$\|u\|_{L^{p}(I,X)}=\Bigg(\int_{I}\|u(t)\|_X^{p}dt\Bigg)^{\frac{1}{p}}$$
      and
      $$\|u\|_{W^{1,p}(I,X)}=\Bigg(\int_{I}\left\|\frac{du}{dt}(t)\right\|_X^{p}+\|u(t)\|_X^{p}dt\Bigg)^{\frac{1}{p}}.$$

\begin{theorem}\label{t:AmannRegularity}
Consider the linear parabolic problem
\begin{equation}\label{e:linear}
    \left\{
   \begin{array}{lr}
     \frac{\partial}{\partial t}u-{\Delta}u+\lambda_{0} u=f \,\,\,\,\,\,\,  in\ \W ,\\
     u(t,x)=0\,\,\,\,\,\,\,\,\mbox{on}\,\,\partial\Omega,\\
     u(0,x)=u_{0}(x) \,\,\,\,\,\,\,\,\mbox{in}\,\,\W,
   \end{array}
   \right.
\end{equation}
where $\lambda_{0}>0$. Given a compact interval $I=[0,T]$, $f\in L^{p}(I,L^{q}(\W))$ and $1<p,q<\infty$,  the solution $u$ to the Problem (\ref{e:linear}) satisfies
\begin{align}\label{inequ:maximalregularity}
\|u\|_{W^{1,q}(I,L^{p}(\W))}+\|u\|_{L^{q}(I,W^{2,p}(\W))}\leq C_{MR}\big(\|u_{0}\|_{W^{s,p}(\W)}+\|f\|_{L^{q}(I,L^{p}(\W))}\big),
\end{align}
where $C_{MR}$ is a positive constant independent of $f$, $u_{0}$ and $I$ and $s>2\left(1-\frac{1}{q}\right)$.
\end{theorem}

\br
In fact, this is a special case of \cite[Theorem \uppercase\expandafter{\romannumeral3}.4.10.7]{Amann}.
We only give this version for our purpose here.
In the original version stated in \cite{Amann}, the first term on the right hand side of (\ref{inequ:maximalregularity}) is in the form of $\|u_{0}\|_{X_{p,q}}$, where the interpolation space $X_{p,q}=(L^{p}(\W),W^{2,p}(\W))_{1-\frac{1}{q},q}$ satisfies $W^{s,p}(\W)\hookrightarrow X_{p,q}$ for $s>2\left(1-\frac{1}{q}\right)$. We refer \cite{LM,Lu,T} once again for details on interpolation spaces.
\er

Now we prove the $H^1$-boundedness of the global solutions.

\begin{lemma}\label{l:H1bounded}
Let $U(t)$ be a global solution to Problem (\ref{e:A14}) with $U(0)\in (H^s)^N$ for $s\in(1,2]$ such that $\lim_{t\to\infty} J(U(t))\geq0$.
Then there is a constant $C>0$ depending only on the $H^2$-norm of $U(0)$ and the initial energy $E_{0}$ such that $\|U(t)\|\leq C$ for any $t\geq0$.
\end{lemma}
\noindent{\bf Proof.}
Denote the interval $I=[t_{0},t_{0}+T]$. Firstly, using the global $L^2$ bound  of $U(t)$ and the same computation of (\ref{inequality:important!!}), we have,
\begin{align}\label{inequality:H1bounded11111111111111}
\sum_{j=1}^{N}\|u_{j}\|^2\leq C(C_{1})\big(1+|\partial_{t}U|_2\big).
\end{align}
Here $C_1$ is the $L^2$ bound of the solution $U(t)$. Using Theorem \ref{t:AmannRegularity} with respect to each equation in Problem (\ref{e:A14}) and putting $p=2$ and $q=\frac{4}{3}$ (Therefore, $s\in\left(\frac{1}{2},2\right]$), we have
$$\|u_{j}\|_{L^{\frac{4}{3}}(I,H^{2}(\W))}\leq C_{MR}\left(\sum_{i=1}^{N}\|U(t_{0})_{i}\|_{H^s}+\Big\|\mu_{j}u_{j}^3+\sum_{i\neq j}\beta_{ij}u^{2}_{i}u_{j}\Big\|_{L^{\frac{4}{3}}(I,L^{2}(\W))}\right).$$
The last term can be estimated as follows:
\begin{align}
\Big\|\mu_{j}u_{j}^3+\sum_{i\neq j}\beta_{ij}u^{2}_{i}u_{j}\Big\|_{L^{\frac{4}{3}}(I,L^{2}(\W))} &\leq C\Bigg(\int_{t_{0}}^{T+t_{0}}\Big(|u_{j}|_{6}^{3}+\sum_{i=1}^{N}|u_{i}^2u_{j}|_{2}\Big)^{\frac{4}{3}}ds\Bigg)^{\frac{3}{4}}\nonumber\\
&\leq C\Bigg(\int_{t_{0}}^{T+t_{0}}\Big(\sum_{i=1}^{N}|u_{i}|_{6}^{3}\Big)^{\frac{4}{3}}ds\Bigg)^{\frac{3}{4}}\nonumber\\
&\leq C\Bigg(\int_{t_{0}}^{T+t_{0}}\Big(\sum_{i=1}^{N}\|u_{i}\|^{4}\Big)ds\Bigg)^{\frac{3}{4}}\nonumber\\
&\leq C\Bigg(\int_{t_{0}}^{T+t_{0}}\Big(\sum_{i=1}^{N}\|u_{i}\|^{2}\Big)^{2}ds\Bigg)^{\frac{3}{4}}.\nonumber
\end{align}
Using (\ref{inequality:H1bounded11111111111111}) and Proposition \ref{p:L2partial_t}, we have
\begin{align}
\Big\|\mu_{j}u_{j}^3+\sum_{i\neq j}\beta_{ij}u^{2}_{i}u_{j}\Big\|_{L^{\frac{4}{3}}(I,L^{2}(\W))} &\leq C(C_{1})\Big(\int_{t_{0}}^{t_{0}+T} \big(1+|\partial_{t}U|_{2}^{2}\big)ds\Big)^{\frac{3}{4}}\nonumber\\
&\leq C(C_{1})(T+E_{0})^{\frac{3}{4}}\leq C(C_{1},E_{0})(T+1)^{\frac{3}{4}}.\nonumber
\end{align}

Now set
\begin{align}\label{constant222222222}
C_{2}=8(NC_{MR})^{2}\left((\|U(0)\|_{(H^s)^N}+1)^2+C(C_{1},E_{0})^{2}(2NC_{MR}+1)^{2}\right)+C(C_{1},E_{0})+1.
\end{align}
Let $T=\big(2NC_{MR}+1\big)^{\frac{4}{3}}$ and $t_{0}=0$. Then we notice that $\|U(0)\|_{(H^s)^N}\leq C_{2}$. And we have
\begin{align}
\Bigg(\int_{0}^{T}\|U\|_{(H^{2})^N}^{\frac{4}{3}}ds\Bigg)^{\frac{3}{4}}&\leq \sum_{j=1}^{N} \Bigg(\int_{0}^{T}\|u_{j}\|_{(H^{2})^N}^{\frac{4}{3}}ds\Bigg)^{\frac{3}{4}}\nonumber\\
&\leq NC_{MR}\Big(\|U(0)\|_{(H^{s})^{N}}+C(C_{1},E_{0})(T+1)^{\frac{3}{4}}\Big).\nonumber
\end{align}
Therefore, there must be a positive number $t'\in(0,T)$ such that
\begin{align}
\|U(t')\|_{(H^{2})^{N}}&\leq NC_{MR}\frac{\|U(0)\|_{(H^s)^N}}{T^{\frac{3}{4}}}+NC(C_{1},E_{0})C_{MR}\Bigg(1+\frac{1}{T}\Bigg)^{\frac{3}{4}}\nonumber\\
&\leq\frac{NC_{MR}\Big(\|U(0)\|_{(H^s)^N}+C(C_{1},E_{0})\Big)}{T^{\frac{3}{4}}}+NC_{MR}C(C_{1},E_{0})\leq C_{2}.\nonumber
\end{align}
We may assume $t'$ is the largest such number in $(0, T]$. With above results, exchanging $t_{0}=0$ into $t_{0}=t'$. Note that we can select $s=2$ for the second and later steps. Via the same method, we can find a largest $t''\in(t',t'+T]$ such that $\|U(t'')\|_{(H^2)^N}\leq C_2$.
Inductively, we can find a sequence of $(t'_{l})_{l}$ such that
\begin{itemize}
  \item $0<t'_{l}-t'_{l-1}\leq T$;
  \item $\lim_{l\to\infty}t'_{l}=\infty$;
  \item $\|U(t'_{l})\|_{(H^2)^N}\leq C_{2}$.
\end{itemize}
The first and the last assertions are obvious. For $\lim_{l\to\infty}t'_{l}=\infty$, we first observe
\begin{align}
\int_{0}^{T}\|U(s)\|_{(H^2)^N}^{\frac{4}{3}}ds&\leq4(NC_{MR})^2\left(\|U(0)\|_{(H^s)^N}^2 +C(C_{1},E_{0})^2\Big((2NC_{MR}+1)^{\frac{4}{3}}+1\Big)\right)\nonumber\\
&\leq C_{2}.\nonumber
\end{align}
This implies that
\begin{align}
C_{2}&\geq\int_{0}^{T}\|U(s)\|_{(H^s)^N}^{\frac{4}{3}}=\int_{\|U\|_{(H^2)^N}<C_{2}}+\int_{\|U\|_{(H^2)^N}\geq C_{2}}\|U(s)\|_{(H^2)^N}^{\frac{4}{3}}ds\nonumber\\
&\geq0+(T-\delta)C_{2}^{\frac{4}{3}}.\nonumber
\end{align}
where $\delta=|\{t\in[0,T]|\|U(t)\|_{(H^2)^N}<C_{2}\}|$. This gives $\delta\geq T-C_{2}^{-\frac{1}{3}}>0$. Therefore for any $l=0,1,\dots$, we have $t'_{l+1}-t'_{l}\geq\delta>0$.

Using the method in Section 2.2 on every interval $[t'_{l},t'_{l}+T]$, we can prove that $\|U(t)\|_{(H^1)^N}\leq C(C_{2})=C(\|U(0)\|_{(H^s)^N})$ for any $t\in [t'_{l},t'_{l}+T]$ and any $l=1,2,\dots$. Therefore, $\|U(t)\|_{(H^1)^N}$ is bounded for $t\geq0$ and the upper bound dependence on $\|U(0)\|_{(H^2)^N}$ continuously.
\begin{flushright}
\qedsymbol
\end{flushright}

We now give two corollaries which will be useful in following paragraph.

\begin{corollary}
For any $U\in\partial\mathcal{A}$, $\|\eta^{t}(U)\|_{(H^1)^N}\leq C$ for any $t\geq0$. Here, the constant $C>0$ depends continuously on the initial data.
\end{corollary}

\br
With the results in Remark \ref{remark:globalexistence}, if we assume $U(0)\in(H^s)^N$, then  we can conclude $\eta^t(U)\in L^{\infty}((0,\infty),(H^s)^N)$ for $s\in\left[1,\frac{3}{2}\right)\cap\left(\frac{3}{2},2\right]$ via the formula of variation of constants.
\er




\subsection{Finer Nodal Properties}

It is well known that for scalar equations along the heat flow the number of changing sign is non-increasing (\cite{An, CMT, Mat}).
For coupled systems it was proved in \cite{DWW} for two equations.
We can prove this is also the case for our system \label{e:Al3} using the arguments in  \cite{CMT} and \cite{DWW},
and we omit details here. But we need a more specific version of this theorem from \cite{CMT}, which is based on the notation of bumps of a radial function.

Let us recall
this from \cite{CMT}.
 The number of changing sign of a continuous radial function $u=u(|x|)$, denoted by $n(u)$, is defined as the largest number $k$ such that there exist a sequence of real number $0<x_{0}<x_{1}<\dots<x_{k}$ such that
      $$u(x_{j})\cdot u(x_{j+1})<0,\,\,\,\,\,\,j=0,\dots,k-1.$$
      We call $n(u)$ the nodal number of the function $u$.
We always assume that the functions we discussed have finite nodal numbers. For a radial function $u$ with $n(u)=k$ and $u(x_{0})>0$, we define its $q$-th bump for $q=1,...,k+1$, by
      \begin{align}
      u_{1}(x)&=\chi_{\{u>0\}}\cdot\chi_{\{|x|<x_{1}\}}\cdot u(x),  \nonumber\\
      u_{q}(x)&=\chi_{\{(-1)^{q-1}u>0\}}\cdot\chi_{\{x_{q-2}<|x|<x_{q}\}}\cdot u(x),\,\,\,\,q=2,\dots,k+1.\nonumber
      \end{align}
      For a radial function $u$ with $n(u)=k$ and $u(x_{0})<0$, we define its $q$-th bump $q=1,...,k+1$ by
      \begin{align}
      u_{1}(x)&=\chi_{\{u<0\}}\cdot\chi_{\{|x|<x_{1}\}}\cdot u(x),  \nonumber\\
      u_{q}(x)&=\chi_{\{(-1)^{q-1}u<0\}}\cdot\chi_{\{x_{q-2}<|x|<x_{q}\}}\cdot u(x),\,\,\,\,q=2,\dots,k+1.\nonumber
      \end{align}
      To avoid confusion, for the $j$-th component $u_{j}$ of $U=(u_{1},\dots,u_{N})$, we denote its $q$-th bump by $u_{j,q}$.

For the solution $U(t)$ to Problem (\ref{e:A14}) with initial value $U\in (H_{r}^{2})^{N}$ for $t\in[0,T(U))$, we denote its $j$-th component by $u(t)_{j}$. By $u(t)_{j,q}$ we denote the $q$-th bump of its $j$-th component.

In this subsection, we always assume that  the initial data $U(0)\in(H^2)^N$. Theorem \ref{t:2} ensure that $\eta^t(U)\in(H^2)^N$ for any $t\geq0$.

We firstly consider the case for $\beta_{ij}\leq0$ for all $i,j=1,\dots,N$ and $i\neq j$.

\begin{proposition}\label{p:L4}
There is a positive number $\rho>0$ such that if $|u_{j,q}|_{4}<\rho$ then $|u_{j,q}(t)|_{4}<\rho$ for $t\geq0$.
\end{proposition}
\noindent{\bf Proof.}
By Theorem \ref{t:2} and the inclusion $H^{2}\subset C(\W)$ for dimension $n=2,3$, then $U(t)$ is continuous in spatial variable. As a consequence, the nodal number of $U(t)$, $n\left(U(t)\right)$, is well-defined. Hence, there exists a small $\varepsilon>0$ such that if
$$(-1)^{q+1}u_{j}(x_{q},0)>0,$$
then
$$(-1)^{q+1}u_{j}(x_{q},t)>0$$
for any $t\in[0,\varepsilon]$. Hence, due to the definition of bump $u_{j,q}$ (c.f. Section 2.1), the differential $\frac{\partial}{\partial t}\int|u_{j,q}|^{4}$ is well-defined. Notice that
\begin{align}
\frac{\partial}{\partial t}\int|u_{j,q}|^{4} &=4\int u_{j,q}^{3}\partial_{t}u_{j,q}=4\int u_{j,q}^{3}\partial_{t}u_{j}\nonumber\\
&=4\int u_{j,q}^{3}\Big(\Delta u_{j}-\lambda_{j} u_{j}+\mu_{j} u_{j}^{3}+\sum_{i\neq j}\beta_{ij} u_{j}u_{i}^{2}\Big)\nonumber\\
&=-3\int |\nabla(u_{j,q}^{2})|^{2}-4\lambda_{j} \int u_{j,q}^{4}+4\mu_{j} \int u_{j,q}^{6}+4\sum_{i\neq j}\beta_{ij} \int u_{j,q}^{4}u_{i}^{2}.\nonumber
\end{align}
Denote $W=u_{j,q}^{2}$. By computing
$$\frac{1}{3}=\frac{\frac{1}{2}}{6}+\frac{1-\frac{1}{2}}{2},$$
 we have from Sobolev embedding of dimensions 2 and 3,
$$|W|_{3}^{3}\leq C\|W\|^{\frac{3}{2}}|W|_{2}^{\frac{3}{2}}.$$
Therefore,
\begin{align}
\frac{\partial}{\partial t}\int |u_{j,p}|^{4}&\leq-C\|W\|^{2}+C\|W\|^{\frac{3}{2}}|W|^{\frac{3}{2}}\nonumber\\
&\leq-C\|W\|^{\frac{3}{2}}|W|_{2}^{\frac{1}{2}}+C\|W\|^{\frac{3}{2}}|W|^{\frac{3}{2}}\nonumber\\
&=-C\|W\|^{\frac{3}{2}}|W|_{2}^{\frac{1}{2}}\big(1-C|W|_{2}\big)\nonumber\\
&=-C\|W\|^{\frac{3}{2}}|W|_{2}^{\frac{1}{2}}\big(1-C|u_{j,q}|_{4}^{2}\big)<0\nonumber
\end{align}
for $|u_{j,q}|_{4}$ small enough.
\begin{flushright}
\qedsymbol
\end{flushright}

If there is a couple $(i_{0},j_{0})$ such that $i_{0}\neq j_{0}$, $i_{0},j_{0}=1,\dots,N$ and $\beta_{i_{0},j_{0}}>0$, the property becomes more
delicate.

\begin{lemma}\label{l:invariant2}
For a solution $U(t)$ with its initial data $U(0)\in(H^2)^N$, then there is a $b=b(\|U(0)\|_{(H^2)^N})>0$ such that if $\beta_{ij}<b$ for all $i,j=1,\dots,N$ and $i\neq j$, the conclusion of the last lemma holds true.
\end{lemma}
\noindent{\bf Proof.}
With the same computation, we have
\begin{align}
\frac{\partial}{\partial t}\int|u_{j,q}|^{4} &=-3\int |\nabla(u_{j,q}^{2})|^{2}-4\lambda_{j} \int u_{j,q}^{4}+4\mu_{j} \int u_{j,q}^{6}+4\sum_{i\neq j}\beta_{ij} \int u_{j,q}^{4}u_{i}^{2}\nonumber\\
&\leq-C_{0}\|W\|^2+4\mu_{j}|W|_{3}^{3}+4\max_{ij}\beta_{ij}\sum_{i=1}^{N}\int W^2 u_{i}^2 ,\nonumber
\end{align}
where $C_{0}=\min\{3,4\lambda_{j}\}$ and $W=u_{j,q}^2$. Now we deal with the last term and have
\begin{align}
4\max_{ij}\beta_{ij}\sum_{i=1}^{N}\int W^2 u_{i}^2 &=4\max_{ij}\beta_{ij}\int W^2 \Big(\sum_{i=1}^{N}u_{i}^2\Big)\leq 4N\max_{ij}\beta_{ij}|W|_{3}^{2}|U|_{6}^{2}\nonumber\\
&\leq 4N\max_{ij}\beta_{ij}S_{3}^{2}S_{6}^{2}\|W\|^{2}\|U\|^{2}\nonumber\\
&\leq 4NbC_{3}(U(0))S_{3}^{2}S_{6}^{2}\|W\|^{2},\nonumber
\end{align}
where $S_{p}$ is the best constant for the inequality $|U|_{p}\leq S_{p}\|U\|$ and $C_{3}(U(0))$ is the upper bound of $\|U(t)\|^{2}$ we computed in the last subsection. If we assume that $b=\frac{\min\{3,4\lambda_{j}\}}{8NS_{3}^{2}S_{6}^{2}C_{3}(U_{0})}>0$, we will have
$$\frac{\partial}{\partial t}\int|u_{j,q}|^{4}\leq-\frac{C_{0}}{2}\|W\|+4\mu_{j}|W|_{3}^{3}.$$
The rest part of the proof is the same with the last lemma.
\begin{flushright}
\qedsymbol
\end{flushright}

\section{Proof of Theorem \ref{t:exist}}

\subsection{A Topological Lemma}

We give the linking structure without assuming any symmetry, which can be used in the proof of Theorem \ref{t:exist}. The strategy of proving is to extend the original setting into a symmetric setting.

\begin{lemma}\label{l:linking}
For $A=[0,+\infty)^{n}$ and a bounded open neighbourhood $\mathcal{O}$ of the origin $0$ in $\mathbb{R}^n$, there is no continuous map $F:\partial\mathcal{O}\cap A\to A$ such that $F(\partial\mathcal{O}\cap A)\subset\partial A\backslash\{0\}$ and the condition $x_{j}=0$ implies that $F_{j}(x_{1}\dots,x_{n})=0$, where $F_{j}(x_{1}\dots,x_{n})$ is the $j$-th component of the vector $F(x_{1}\dots,x_{n})$.
\end{lemma}
\noindent{\bf Proof.}
We argue it by contradiction. Suppose there is a continuous mapping $F:\partial\mathcal{O}\cap A\to A$ such that $F(\partial\mathcal{O}\cap A)\subset\partial A\backslash\{0\}$. Inspired by \cite{LM1111111}, we begin the proof by extending the setting to a symmetric version and obtain the contradiction via a genus argument.

Firstly, let us define another open neighbourhood $\mathcal{O}^*$ of the origin $0\in\mathbb{R}^n$ by reflection with respect to each component of the coordinates, i.e.
$$\mathcal{O}^* =\{x=(x_{1},\dots,x_{n})|(|x_{1}|,\dots,|x_{n}|)\in \mathcal{O}\}.$$
It is easy to see that the open set $\mathcal{O}^*$ is antipodal symmetric. Then the following inclusion holds true:
\begin{align}\label{inclusion:boundary}
\partial\mathcal{O}^{*}\cap A\subset \partial\mathcal{O}\cap A.
\end{align}
Indeed, we observe that $\partial\mathcal{O}^{*}\cap int(A)= \partial\mathcal{O}\cap int(A)$. So we only need to show that $\partial\mathcal{O}^{*}\cap \partial A\subset \partial\mathcal{O}\cap\partial A$. For any $x\in \partial\mathcal{O}^{*}\cap \partial A$, for any $r>0$, $B_{\mathbb{R}^n}(x,r)\cap\mathcal{O}^{*}\neq\emptyset$. Due to the construction of $\mathcal{O}^*$, the last intersection gives $B_{\mathbb{R}^n}(x,r)\cap\mathcal{O}\neq\emptyset$, which implies that $x\in \partial\mathcal{O}\cap \partial A$. Now we restrict the mapping $F$ to the set $\partial\mathcal{O}^{*}\cap A$ and extend it to the whole $\partial\mathcal{O}^{*}$. Define the mapping $\tilde{F}:\partial\mathcal{O}^{*}\to X$ by
$$\tilde{F}(x_{1},\dots,x_{n})=\big(sgn(x_{1})F_{1}(|x_{1}|,\dots,|x_{n}|),\dots,sgn(x_{n})F_{n}(|x_{1}|,\dots,|x_{n}|)\big),$$
where $X=\{x=(x_{1},\dots,x_{n})|\prod_{j=1,\dots,n}x_{j}=0\}$. Then we can claim
\begin{itemize}
  \item $\tilde{F}$ is an odd extension of $F$;
  \item $\tilde{F}(\partial\mathcal{O}^{*})\subset X\backslash \{0\}$.
\end{itemize}
The first assertion is easy. Now we check the second one. We only need to verify that for any $x\in\partial\mathcal{O}^{*}$, $\tilde{F}(x)\neq0$. Otherwise, if for any $j=1,\dots,n$, $sgn(x_{j})F_{j}(|x_{1}|,\dots,|x_{n}|)=0$. Since $x\in\partial\mathcal{O}^{*}$, there are some integers, say $1,\dots,s$, with $x_{1}=\dots,x_{s}=0$ and other integers, say $s+1,\dots,n$, with $x_{s+1}\neq0$, $\dots$, $x_{n}\neq0$. Therefore, $F_{j}(|x_{1}|,\dots,|x_{n}|)\equiv0$ for any $j=1,\dots,n$. This is a contradiction with $(|x_{1}|,\dots,|x_{n}|)\in\partial\mathcal{O}\cap A$ and $F(\partial\mathcal{O}\cap A)\subset\partial A\backslash\{0\}$.

Now we will have a contradiction via the genus with respect to the symmetry of antipodal.
We denote the genus generated by the antipodal symmetry by $\gamma'$.
 On one hand, we have $n=\gamma'(\partial\mathcal{O}^{*})\leq\gamma'(\tilde{F}(\partial\mathcal{O}^{*}))$ due to Borsuk's theorem for the symmetry with respect to antipodal. On the other hand, notice that $X\cap X'=\{0\}$, where $X'=\{x=(x_{1},\dots,x_{n})|x_{1}=\dots=x_{n}\}$. Then we can construct an odd homotopy $G$ such that
$$X\backslash\{0\}\overset{G}{\simeq}\mathbb{S}^{n-2}.$$
This implies that $\gamma'(\tilde{F}(\partial\mathcal{O}^{*}))\leq n-1$. This is a contradiction.
\begin{flushright}
\qedsymbol
\end{flushright}
We remark that some of the computations in the above proof were used in \cite{LW1,LW} and will be also used in the next section.

\subsection{Proof of Theorem \ref{t:exist}}

We prove this theorem via the concept of invariant sets of a descending flow and we will use the parabolic flow as a mean of descending flow.
Recall we fix $N$ non-negative integers $P_{1},\dots,P_{N}$ which are the prescribed componentwise nodal numbers.
We first introduce some auxiliary functions.
\begin{itemize}
  \item Firstly, for the radial domain $\W$, we cut it into $N$ radial sub-domains $\W_{j}$ for $j=1,\dots,N$ with $\overline{\W}=\cup_{j=1}^{N}\overline{\W_{j}}$;
  \item For any fixed $j=1,\dots,N$, we cut the domain $\W_{j}$ into $P_{j}+1$ sub-domains $\W_{j,q}$ with $q=1,\dots,P_{j}+1$ with $\overline{\W_{j}}=\cup_{q=1}^{P_{j}+1}\overline{\W_{j,q}}$;
  \item For any $j=1,\dots,N$ and $q=1,\dots,P_{j}+1$, we define a smooth non-zero radial function with compact support $w_{j,q}:\W_{j,q}\to[0,+\infty)$.
\end{itemize}
Without loss of generality, we can assume that $|w_{j,q}|_{4}\equiv 1$ for any $j=1,\dots,N$ and $q=1,\dots,P_{j}+1$. We define the following set $S$ by
\begin{align}
S=\Bigg\{\Bigg(\sum_{q=1}^{P_{1}+1}(-1)^{q+1}\alpha_{1,q}w_{1,q}(x),&\dots,\sum_{q=1}^{P_{N}+1}(-1)^{q+1}\alpha_{N,q}w_{N,q}(x)\Bigg)\Bigg|\nonumber\\ &\alpha_{j,q}\geq\frac{\varepsilon}{100}\,\,\mbox{for}\,\,q=1,\dots,P_{j}+1\,\,and\,\,j=1,\dots,N\Bigg\},\nonumber
\end{align}
which is a closed cone in the real Euclidean space of dimension $\sum_{j=1}^{N}P_{j}+N$. And there is an isomorphism
\begin{align}
i:S&\to[0,+\infty)^{\sum_{j=1}^{N}P_{j}+N}\nonumber\\
\left(\sum_{q=1}^{P_{j}+1}\alpha_{j,q}w_{j,q}\right)_{j}&\mapsto\Big(\alpha_{j,q}-\frac{\varepsilon}{100}\Big)_{j,q}.\nonumber
\end{align}
We denote by $Y$ the $\sum_{j=1}^{N}P_{j}+N$ dimensional Euclidean space spanned by $S$ with respect to the linearity in $i(S)$. It is easy to see that $\mathcal{A}\cap Y$ is also an open neighbourhood of the origin in $Y$ and $\mathcal{A}\cap Y$ is bounded. Notice that $S\cap\partial\mathcal{A}$ is a compact set and can be embedded into a finite dimensional Euclidean space, where all the norms are equivalence. If there is at least one $\beta_{ij}$ is positive, due to the limit of Lemma \ref{l:invariant2}, we need to find an upper bound $b>0$ depending on $\sup_{U\in S\cap\partial\mathcal{A}}\sup_{t\geq0}\|\eta^{t}(U)\|^2\in(0,\infty)$. If all the $\beta_{ij}$'s are non-positive, the limit on the upper bound is no longer necessary (Lemma \ref{p:L4}).

Now we will locate the portions of the boundary $\partial \mathcal A$ in which along the flow lines the number of nodal domains can be controlled. For this purpose, as done for scalar equations in \cite{CMT}, we introduce the following notations:
\begin{itemize}
  \item $D_{j,k}=\{U=(u_{1},\dots,u_{N})\in(H_{r}^{2})^{N}|n(u_{j})=k\}$ and $D=\cap_{j=1}^{N}D_{j,P_{j}}$, where $P_{1},\dots,P_{N}$ are given in the theorem;
  \item $E_{j,q}^{\varepsilon}=\{U=(u_{1},\dots,u_{N})\in D||u_{j,q}|_{4}<\varepsilon\}$ for $q=1,\dots,P_{j}$ and $j=1,\dots,N$;
  \item denote
  \begin{align}
  H=\Bigg\{U=(u_{1},\dots,u_{N})\in(H_{r}^{2})^{N}|& n(u_{j})\leq P_{j}\,\,\mbox{for}\, j=1,\dots,N\nonumber\\
   &\mbox{and}\, \sum_{j=1}^N n(u_j) < \sum_{j=1}^N P_j\Bigg\}\nonumber
  \end{align}
  \item denote
  \begin{align}
  F_{\varepsilon}=\cup_{j=1}^{N}\cup_{q=1}^{P_{j}}E_{j,q}^{\varepsilon}\cup H;\nonumber
  \end{align}
  \item the complete invariant set of the set $E_{j,q}^{\varepsilon}$ is defined as
  $$C(E_{j,q}^{\varepsilon})=\big\{U\in(H^{2})^{N}\big|\exists t_{0}\geq0\,\,s.t.\,\,\eta^{t_{0}}(U)\in E_{j,q}^{\varepsilon}\big\}$$
  for $q=1,\dots,P_{j}+1$ and $j=1,\dots,N$. Therefore, we can denote
  $$A^{\varepsilon}_{j,q}=C(E_{j,q}^{\varepsilon})\cap\partial\mathcal{A}\cap D$$
  for $q=1,\dots,P_{j}+1$ and $j=1,\dots,N$.
\end{itemize}
Due to the invariance property proved in the last section we will also define an arriving time
\begin{itemize}
  \item for any $U\in D\cap\partial\mathcal{A}$, denote $T^{*}(U)=\inf\{t\geq0|\eta^{t}(U)\in F_{\frac{\varepsilon}{2}}\}$.
\end{itemize}

\vskip .1in

We note that the set $D_{j,k}$ consists of the vector-valued functions whose $j$-th components have exactly $k$ sign-changing number. For any function $(u_{1},\dots,u_{N})$ in the set $D$, let $u_{j}$ be an arbitrary component. The set $E_{j,q}^{\varepsilon}$ contains the functions whose $q$-th bump of $j$-th component has a small $L^{4}$ norm. This set is invariant due to Proposition \ref{p:L4} when $n(\eta^{t}(\cdot)_{j})$ dose not change. As to the set $F_{\varepsilon}$, an element in $F_{\varepsilon}$ either has a small bump or has a component with a sign-changing number less than what we prescribed. It should be noted that this set is what we want to remove. $T^{*}$ is the time when the flow line arrives in the set $F_{\frac{\varepsilon}{2}}$. Using the computation in Proposition \ref{p:L4}, any flow line which flows into $F_{\varepsilon}$ in a finite time will flow into $F_{\frac{\varepsilon}{2}}$ eventually.
We will fix small $\varepsilon$ so that the invariance property holds.
We now have the continuity of the arriving time $T^* (U)$.
\begin{lemma}
The arriving time $T^* (U)$ is continuous.
\end{lemma}
\noindent{\bf Proof.}
Due to the detailed computation in the following paragraph, we restrict ourselves to the case $n(u_{j})\leq P_{j}$ for $j=1,\dots,N$. Here $u_{j}$ denotes the $j$-th component of the vector-valued function $U$. In fact, we only need the case $n(u_{j})\equiv P_{j}$ for any $j=1,\dots,N$ in this section. Nonetheless, we still prove the general case for the sake of the next section.

First we consider the case $\sum_{j=1}^{N}n(u_{j})<\sum_{j=1}^{N}P_{j}$, so we have $T^{*}(U)=0$. We prove the continuity at $U$ by contradiction.
For a sequence $U_{n}\to U$ in $(H^{1})^{N}$ with $U_{n}\in D$ for any $n=1,2,\dots$, suppose there is $t_{0}>0$ such that $T^{*}(U_{n})\geq 2t_{0}$ for large $n$. We select and fix a $t\in(0,t_{0})$. On one hand, we have
$$\eta^{t}(U_{n})\to \eta^{t}(U)\qquad in\,\,(H^{1})^{N},$$
due to Theorem \ref{t:2}. On the other hand, by the definition of the arriving time $T^{*}$ and the non-increasing property of nodal number along the flow line, using $U_{n}\in D$ we have
\begin{itemize}
  \item $n\big((U_{n})_{j}\big)=P_{j}$ for $j=1,\dots,N$;
  \item $\big|\big(\eta^{t}(U_{n})\big)_{j,q}\big|_{4}\geq\frac{\varepsilon}{2}$ for $j=1,\dots,N$ and $q=1,\dots,P_{j}$;
  \item $n\big(\big(\eta^{t}(U)\big)_{j}\big)\leq P_{j}$ for $j=1,\dots,N$ and $q=1,\dots,P_{j}$ and at least one of the $\leq$'s holds strictly.
\end{itemize}
Here, $(\eta^{t}(W))_{j}$ and $(\eta^{t}(W))_{j,q}$ are the $j$-th component and the $q$-th bump of the $j$-th component of $\eta^{t}(W)$. Now we show that these assertions lead us to a contradiction.

Since $\eta^{t}(U_{n})\to\eta^{t}(U)$ in $(H^{1})^{N}$, we can select a large $n_{0}>0$ such that $\big|\eta^{t}(U_{n_{0}})-\eta^{t}(U)\big|_{4}\leq\frac{\varepsilon}{4}$. In the following, we will argue it in terms of components. Let us consider $\eta^{t}(U)_{1}$ and $\eta^{t}(U_{n_{0}})_{1}$ for the sake of simplicity, where $\eta^{t}(U)_{1}$ and $\eta^{t}(U_{n_{0}})_{1}$ are the first components of $\eta^{t}(U)$ and $\eta^{t}(U_{n_{0}})$ respectively. Let us assume that $n(\eta^{t}(U)_{1})<n(\eta^{t}(U_{n_{0}})_{1})=P_{1}$ without loss of generality. Due to the definition of the sign-changing number, we can find a sequence of numbers $x_{q-1}\in \mbox{supp} \,\eta^{t}(U_{n_{0}})_{1,q}$ for $q=1,\dots,P_{1}+1$, where $\eta^{t}(U_{n_{0}})_{1,q}$ is the $q$-th bump, such that
$$\eta^{t}(U_{n_{0}})_{1}(x_{q})\cdot\eta^{t}(U_{0})_{1}(x_{q+1})<0$$
for $q=0,\dots,P_{1}$. Using the facts $\big|\eta^{t}(U_{n_{0}})_{1,q}\big|_{4}\geq\frac{\varepsilon}{2}$ for $q=1,\dots,P_{1}$, and $\big|\eta^{t}(U_{n_{0}})-\eta^{t}(U)\big|_{4}\leq\frac{\varepsilon}{4}$, we claim there must be $x'_{q-1}\in \mbox{supp} \, \eta^{t}(U_{n_{0}})_{1,q}$ such that
\begin{align}\label{ineq:NODALS}
\eta^{t}(U)_{1}(x'_{q})\cdot\eta^{t}(U)_{1}(x'_{q+1})<0
\end{align}
for $q=0,\dots,P_{1}$. Otherwise, if there is a $q_{0}=1,\dots,P_{1}$ such that
\begin{itemize}
  \item $\eta^{t}(U_{n_{0}})\geq0$;
  \item $\eta^{t}(U)\leq0$ on $\mbox{supp}\eta^{t}(U_{n_{0}})_{1,q_{0}}$;
  \item $|\eta^{t}(u)_{1,q_{0}}|_{4}\geq\frac{\varepsilon}{2}$;
\end{itemize}
then we have
\begin{align}
\frac{\varepsilon}{4}&\geq|\eta^{t}(U_{n_0})-\eta^{t}(U)|_{4}\geq |\eta^{t}(U_{n_0})_{1,q_{0}}-\eta^{t}(U)\cdot\chi_{\mbox{supp}\eta^{t}(U_{n_0})_{1,q_{0}}}|_{4}\nonumber\\
&\geq|\eta^{t}(U_{n_0})_{1,q_{0}}|_4\geq\frac{\varepsilon}{2},\nonumber
\end{align}
which is a contradiction. Here, the function $\chi_{A}$ is the characteristic function of the set $A$. Therefore, (\ref{ineq:NODALS}) holds, i.e., $P_{1}=n(\eta^{t}(U)_{1})$, but this is a contradiction with $n(\eta^{t}(U)_{1})<P_{1}$. The proof of the first case is complete.

Next we consider the case $n(U_{j})= P_{j}$ for $j=1,\dots,N$. For a sequence $U_{n}\to U$ in $(H^{1})^{N}$, we only check the lower limit
$$T^{*}(U)\leq\varliminf_{n\to\infty}T^{*}(U_{n}).$$
The upper one
$T^{*}(U)\geq\varlimsup_{n\to\infty}T^{*}(U_{n})$
can be proved in the same way.

We argue by contradiction again. Suppose, up to a subsequence, we have
$$s:=\lim_{n\to\infty}T^{*}(U_{n})<T^{*}(U)\leq T(U)=\infty.$$
Then we can find a $t\in(s,T^{*}(U))$. Set $V_{n}=\eta^{t}(U_{n})$ and $V=\eta^{t}(U)$. Then, due to Theorem \ref{t:2}, we have $V_{n}\to V$ in $(H^{1})^{N}$.
Since $t<T^{*}(U)$, we have $V=\eta^{t}(U)\in D\cap\partial\mathcal{A}\backslash F_{\frac{\epsilon}{2}}$. Then $|V_{j,q}|_{4}\geq \frac{\varepsilon}{2}$ for some $j=1,\dots,N$, $q=1,\dots,P_{j}$.
By using $T^{*}(U_{n})\to s<t$, we have $V_{n}=\eta^{t}(U_{n})\in F_{\frac{\epsilon}{2}}\cap\partial\mathcal{A}$, which implies $|(V_{n})_{j,q}|_{4}\leq\frac{\varepsilon}{2}$.
Combining these with the fact that $V_{n}\to V$ in $(H^{1})^{N}$, we conclude that $|V_{j,q}|_{4}=\frac{\varepsilon}{2}$. Hence, we obtain $|\eta^{t}(U)_{j,q}|_{4}=|\eta^{T^{*}(U)}(U)_{j,q}|_{4}=\frac{\varepsilon}{2}$ with $t<T^{*}(U)$.
Using Proposition \ref{p:L4}, for any $\theta\in(t,T^{*}(U))$, $|\eta^{\theta}(U)_{j,q}|_{4}\equiv\frac{\varepsilon}{2}$, which is a contradiction with the very proposition itself.

\begin{flushright}
\qedsymbol
\end{flushright}

Finally to prove Theorem  \ref{t:exist}, we only need to show that $$A:=\partial\mathcal{A}\cap D\backslash\big(\cup_{j=1}^{N}\cup_{q=1}^{P_{j}+1}A_{j,q}^{\varepsilon}\big)\neq\emptyset.$$
The rest part of the proof requires a lower bound of the energy functional $J$ on the set $A$ and the fact that the energy functional satisfies the (PS) condition. The second part is obvious and the first part is given by
$$A\subset\partial\mathcal{A}$$
and
$$0\leq\inf_{\partial\mathcal{A}}J\leq\inf_{A}J.$$

Now we verify that $\partial\mathcal{A}\cap D\backslash\big(\cup_{j=1}^{N}\cup_{q=1}^{P_{j}+1}A_{j,q}^{\varepsilon}\big)\neq\emptyset$.

This proof relies heavily on a technique used in \cite{LW1,LW}. We will use the subset $S\subset D$ constructed at the beginning of this subsection and prove the theorem
by proving $\partial\mathcal{A}\cap S\backslash\big(\cup_{j=1}^{N}\cup_{q=1}^{P_{j}+1}A_{j,q}^{\varepsilon}\big)\neq\emptyset$.

Now we argue by contradiction, i.e.,
we assume that $\partial\mathcal{A}\cap S\subset\cup_{j=1}^{N}\cup_{q=1}^{P_{j}+1}A_{j,q}^{\varepsilon}$.
We use $\partial_{Y}S$ to denote the boundary of $S$ with respect to the space of $Y$.
Define a continuous cut-off function $\phi:[0,\infty)\to[0,1]$:
\begin{equation}\label{cutoff}
\phi(x)=\left\{
\begin{aligned}
1 & \qquad & s\geq\varepsilon, \\
0 & \qquad & s\leq\frac{\varepsilon}{2}, \\
\frac{2s}{\varepsilon}-1 & \qquad & s\in\Big(\frac{\varepsilon}{2},\varepsilon\Big).
\end{aligned}
\right.
\end{equation}
Let us define the mapping $h:\partial\mathcal{A}\cap S\to\partial_{Y}S$ by
$$
h\big(U\big)=\Bigg(\sum_{q=1}^{P_{1}+1}\Big(\phi(|\eta^{T^{*}(U)}(U)_{1,q}|_{4})+\frac{\varepsilon}{100}\Big)w_{1,q},\dots, \sum_{q=1}^{P_{N}+1}\Big(\phi(|\eta^{T^{*}(U)}(U)_{N,q}|_{4})+\frac{\varepsilon}{100}\Big)w_{N,q}\Bigg),$$
where $U=(U_{1},\dots,U_{N})\in\partial\mathcal{A}\cap S$.

We are here in position to use Lemma \ref{l:linking}. To do this, we only need to check that $h(\partial_{Y}(\mathcal{A}\cap Y)\cap S)\subset \partial_{Y}S\backslash\{\theta\}$. Firstly, we notice that $\partial_{Y}(\mathcal{A}\cap Y)\subset\mathcal{A}\cap Y$.
Then we claim that for any $U\in\partial\mathcal{A}\cap S$, there is a $j=1,\dots,N$ and a $q=1,\dots,P_{j}+1$ such that $\phi(|\eta^{T^{*}(U)}(U)_{j,q}|_{4})>0$.
If the claim is not true, for any $j=1,\dots,N$ and $q=1,\dots,P_{j}+1$, we have $\int|U_{j,q}|^{4}<\varepsilon^4$ as $t$ goes large. Multiplying $u_{j}$ on the both sides of the $j$-th equation of Problem (\ref{e:A14}) and integrating and summing up with respect to $j$, we have
\begin{align}
\frac{1}{2}\frac{\partial}{\partial t}|U(t)|^{2}_{2}+\|U(t)\|^{2}\leq\sum_{j=1}^{N}\Bigg(\mu_{j}\int|u_{j}(t)|^4 +\sum_{i\neq j}\beta_{ij}\int u_{i}(t)^{2}u_{j}(t)^{2}\Bigg)<C\varepsilon\nonumber
\end{align}
 for $t$ large. Using the openness of $\mathcal{A}$ in $(H^1)^N$ and the invariance of $\partial\mathcal{A}$, $\|U(t)\|\geq C>0$ uniformly for $t\geq0$. This implies that $\frac{\partial}{\partial t}|U(t)|_{2}^{2}\leq-C$ for $t>0$,  a contradiction. Therefore, the mapping $h$ satisfies the condition for $F$ in Lemma \ref{l:linking}, so we have a contradiction on the existence of the mapping $h$. The proof is complete.
\begin{flushright}
\qedsymbol
\end{flushright}

\section{Proof of Theorem \ref{t:main}}

\subsection{The Idea of the Proof}

We are now in position to prove Theorem \ref{t:main}.
We outline our approach briefly first. Using the flow invariance property we reduce the variational problem
to one defined on a subset on the boundary of the stable set of the origin where the nodal number of the functions is
controlled by the componentwise-prescribed nodes.
In order to establish multiple nodal solutions having the same componentwise nodal number we will make use of
the symmetry property imposed in Theorem \ref{t:main}.

More precisely, our problem possesses a $\mathbb Z_p$ symmetry under a cyclic permutation $\sigma:(u_{1},\dots,u_{N})\mapsto(\sigma_{1}(u_{1}),\dots,\sigma_{N}(u_{N}))$
  in $(H^{1}_{0}(\W))^N$ defined by
      \begin{itemize}
        \item $\sigma_{i}(u_{i})=u_{i+1}$ for $i\neq pb$ for $b=1,\dots,B$,
        \item $\sigma_{pb}(u_{pb})=u_{p(b-1)+1}$ for $b=1,\dots,B$.
      \end{itemize}
In other words, we define the permutation $\sigma$ as
      \begin{align}
      \sigma(u_{1},u_{2},\dots,u_{p};&\dots\dots;u_{N-p+1},u_{N-p+2},\dots,u_{N})\nonumber\\
      &=(u_{2},\dots,u_{p},u_{1};\dots\dots;u_{N-p+2},\dots,u_{N},u_{N-p+1}).\nonumber
      \end{align}
      It is easy to see that this can be regarded as a $\mathbb{Z}_{p}$ cyclic group action, our variational functional $J$ is invariant under this action.

  We will use a $\mathbb Z_p$ group action index (or genus), which is from \cite{TW} (see also related works in \cite{Wang, WW}). We summarize some basic property of the index. Let $E$ be a Banach space on which there is an $\mathbb Z_p$ action generated by $\sigma$.
  Let $F_\sigma= \{U\in E\;|\; \sigma U = U \}$ be the set of fixed points of the $\sigma$ action. For a $\sigma$-symmetric compact set $A\subset E\backslash F_\sigma $, the index $\gamma(A)$ is defined as the smallest $m\in\mathbb{N}$ such that there exists a mapping $h:A\to\mathbb{C}^{m}\backslash\{0\}$ with
      $$h(\sigma U)=e^{\frac{2\pi i}{p}}h(U).$$
      If there is no such mapping, we set $\gamma(A)=\infty$. We need the following properties of the Index $\gamma$.

\begin{itemize}
        \item If $A\subset B$, then $\gamma(A)\leq\gamma(B)$;
        \item $\gamma(A\cup B)\leq\gamma(A)+\gamma(B)$;
        \item if $g:A\to E\backslash F_{\sigma}$ is continuous and satisfies
        $g(\sigma(u))=\sigma g(u)$ for all $ u\in A$,
        then $$\gamma(A)\leq\gamma(\overline{g(A)});$$
        \item if $\gamma(A)>1$, then $A$ is an infinite set;
        \item if $A$ is compact and $\gamma(A)<\infty$, then there exist an open $\sigma$-invariant neighbourhood $\mathcal{N}$ of $A$
        such that $\gamma(A)=\gamma(\overline{\mathcal{N}})$;
        \item if $S$ is the boundary of a bounded neighbourhood of the origin in a $m$-dimensional complex linear space such that
       $e^{\frac{2\pi i}{p}}U \in S $ for any $U\in S$, and $\Psi:S\to E\backslash F_\sigma $ is continuous and satisfies for any $U\in S$,
       $\Psi(e^{\frac{2\pi i}{p}}U)=\sigma(\Psi(U))$, then $\gamma(\Psi(S))\geq m$;
\item let $A$ be a closed set such that $\mathbb Z_p(A)\subset E\setminus F_\sigma$, and $\cap_{i=0}^{p-1} \sigma^i(A) =\emptyset $.
Then $\gamma(\mathbb Z_p(A))\leq p-1$.
   \end{itemize}
One can find the proofs of all these properties but the last one in
for example \cite{TW}.
 We sketch the proof of the last dot here.

Note that it is clear that if $A\cap\sigma(A)=\emptyset$ then $\gamma(\mathbb{Z}_{p}(A))=1$, which is analogous to the normalization of genus (c.f. \cite[Proposition 7.5]{R}). We now decompose $\mathbb{Z}_{p}(A)$ into $p-1$ sets of genus equal to $1$. For small $\varepsilon>0$, define$A_{1}:=\cap_{i=0}^{p-2}\sigma^{i}(A)$, $A_{2}:=\cap_{i=0}^{p-3}\sigma^{i}(A)\backslash\left(\cap_{i=0}^{p-2}\sigma^{i}(A)\right)_{\varepsilon}$, ..., $A_{p-1}:=A\backslash\left(A\cap\sigma(A)\right)_{\varepsilon}$. Using the construction of $A_{i}$'s and the condition $\cap_{i=0}^{p-1} \sigma^i(A) =\emptyset$, we can verify that $A_{i}\cap\sigma(A_{i})=\emptyset$ for any $i=1,\dots,p-1$. The conclusion follows from $\mathbb{Z}_{p}(A)\subset\cup_{i=1}^{p-1}\overline{\left(\mathbb{Z}_{p}(A_{i})\right)_{\varepsilon}}$.


Using this $\mathbb Z_p$-cyclical permutation symmetry and the genus type index generated as above, we construct multiple nodal solutions with a given componentwise prescribed node by a minimax type argument.
With the aid of the flow invariance, the central part of the proof is to construct a certain set of vector-valued functions which has infinite $\mathbb Z_p$ genus and in which the flow line always possesses prescribed number of nodal domains. Then by a minimax construction in variational methods (c.f. \cite{Ch0, R}), we will have a sequence of critical levels and therefore a sequence of solutions to Problem (\ref{e:A111}).
For the construction of sets with large $\mathbb Z_p$-index, we will use a variant of the construction in \cite[Proposition 4.2]{TW} where only positive solutions were considered, by making sets of sign-changing functions with $\sigma$ symmetry property.
For computations of the $\mathbb Z_p$ index, we will adapt some ideas from \cite{CMT, LW1, LW, TW} incorporating the invariance of nodal domains and the $\mathbb Z_p$ symmetry.


\subsection{Invariant Sets and Other Constructions}

We need a symmetric version of the settings in Section 3. We begin with constructing sets of vector-valued functions with componentwise-prescribed number of nodal domains and with arbitrarily large genus.

Recall that $p$ is a prime factor of $N$ and $B$ is such that $N=pB$ and that $P_{1},\dots,P_{B}$ are $B$ non-negative integers and fixed in the proof.
For any given positive integer $K$, we will construct a subset having $\mathbb Z_p$ genus not less than $K$ that consists of vector-valued functions $U=(u_{1},\dots,u_{N})$ such that for $b=1,\dots,B$, $n(u_{pb-p+i})=P_{b}$ for $i=1,...,p$, and satisfies other dynamic property.

Firstly, we divide the domain $\W$ into $B$ radial parts and denote them by $\W_{b}$ for $b=1,...,B$ so $\overline{\cup_{b=1}^{B}\W_{b}}=\overline{\W}$. For a fixed integer $b=1\dots,B$, we divide $\W_{b}$ into $P_{b}+1$ radial sub-domains $\W_{b,q}$ for $q=1,\dots,P_{b}+1$ so $\overline{\cup_{q=1}^{P_{b}+1}\W_{b,q}}=\overline{\W_{b}}$. For each sub-domain $\W_{b,q}$, divide it into $K$ radial sub-domains $\W_{b,q,k}$ for $k=1,\dots,K$ so $\overline{\cup_{k=1}^{K}\W_{b,q,k}}=\overline{\W_{b,q}}$. Denote $\mathcal{O}_{b,q}=\mathbb{S}^{1}\times\W_{b,q}$ and $\mathcal{O}_{b,q,k}=\mathbb{S}^{1}\times\W_{b,q,k}$. For $b=1\dots,B$, $q=1,\dots,P_{b}$ and $k=1,\dots,K$, we define functions for $(t,x) \in \mathcal{O}_{b,q,k}=\mathbb{S}^{1}\times\W_{b,q,k}$ as follows:
\begin{itemize}
  \item $w_{b,q,k}(t,x)=w_{b,q,k}(t,|x|)=w_{b,q,k}(t,r):\mathcal{O}_{b,q,k}\to\mathbb{R}$ of class $C^{4}$ and of compact support in $\mathcal{O}_{b,q,k}$;
  \item $w_{b,q,k}\geq0$ and $w_{b,q,k}(t,\cdot)\not\equiv0$ for any $t\in\mathbb{S}^1$ ;
  \item $\mbox{supp}\, w_{b,q,k}(t,\cdot)\cap \mbox{supp}\, w_{b,q,k}\big(\frac{2\pi }{p}+t,\cdot\big)=\emptyset$ for any $t\in\mathbb{S}^1$.
\end{itemize}

Now we have a few words about the notation for clarity. We point out that the subscript "$b$" denotes the number of the blocks of components with each block having $p$ components so it is invariant under the $\mathbb{Z}_{p}$-permutation of components. The subscript "$q$" denotes the number of nodal domains. And the subscript "$k$" is for the factor $K$ of the dimension of the simplex.
 To give a simplex in Sobolev space $(H^{1})^{N}$ involving vector-valued functions, we start by its componentwise construction. In order to use the $\mathbb Z_p$ index $\gamma$, we need to consider the complex Euclidean space $\mathbb{C}^{K\big(\sum_{b=1}^{B}P_{b}+B\big)}$. For any vector $z=(z_{b,q,k})$ we decompose them in the polar-coordinates with respect to the components. This leads to $z_{b,q,k}=e^{i\theta_{b,q,k}}\alpha_{b,q,k}$ with $\alpha_{b,q,k}$'s are nonnegative real numbers and $\theta_{b,q,k}\in[0,2\pi)$ for any $b=1\dots,B$, $q=1,\dots,p$ and $k=1,\dots,K$.  For $b=1,...,B$ fixed we define
\begin{align}
V_{b}(t,z_{b})(r)=\sum_{q=1}^{P_{b}+1}(-1)^{q+1}\sum_{k=1}^{K}\alpha_{b,q,k}w_{b,q,k}(t+{\theta_{b,q,k}},r)\nonumber
\end{align}
where the vector $z_{b}=\{(z_{b,q,k})\,|\, q=1,\dots,P_{b}+1; k=1,\dots,K\}$. Then we can define a mapping $$\psi:\mathbb{C}^{K\big(\sum_{b=1}^{B}P_{b}+B\big)}\to(H_{r}^{2})^{N}$$ by
\begin{align}
\psi(z)=\Bigg(V_{1}(0,z_{1}),&V_{1}\Big(\frac{2\pi}{p},z_{1}\Big),\dots,V_{1}\Big(\frac{2\pi (p-1)}{p},z_{1}\Big),\nonumber\\
&\dots\qquad\dots\nonumber\\
&V_{B}(0,z_{B}),V_{B}\Big(\frac{2\pi}{p},z_{B}\Big),\dots,V_{B}\Big(\frac{2\pi (p-1)}{p},z_{B}\Big)\Bigg).\nonumber
\end{align}
We note that
\begin{align}
V_{b}(t,e^{\frac{2\pi i}{p}}z_{b})(r)&=\sum_{q=1}^{P_{b}+1}(-1)^{q+1}\sum_{k=1}^{K}\alpha_{b,q,k}w_{b,q,k}\Big(t+\theta_{b,q,k}+\frac{2\pi}{p},r\Big)\nonumber\\
&=V_{b}\Big(\frac{2\pi}{p}+t,z_{b}\Big)(r),\nonumber
\end{align}
which implies that $\psi\big(e^{\frac{2\pi i}{p}}z\big)=\sigma\psi(z)$. Here recall $\sigma$ is the $\mathbb Z_p$ cyclic permutation.

As the settings and notations used in Section 3, we introduce the following notations:
\begin{itemize}
  \item $D_{j,k}=\{U=(u_{1},\dots,u_{N})\in(H_{r}^{2})^{N}|n(u_{j})=k\}$ and $$D=\cap_{b=1}^{B}\cap_{i=1}^{p}D_{pb-p+i,P_{b}};$$
  \item For $\varepsilon >0$, $E_{j,q}^{\varepsilon}=\{U=(u_{1},\dots,u_{N})\in D||u_{j,q}|_{4}<\varepsilon\}$;
  \item denote
  \begin{align}
  H=\Bigg\{U=(u_{1},\dots,u_{N})\in(H_{r}^{2})^{N}|& n(u_{bp-p+i})\leq P_{b},\,\mbox{for}\, i=1,\,\dots,p,\,\,b=1,\dots,B,\;\nonumber\\
   &\mbox{and}\, \sum_{j=1}^N n(u_j) < p\sum_{b=1}^B P_b\Bigg\}\nonumber
  \end{align}
  and
  \begin{align}
  F_{\varepsilon}=&\cup_{b=1}^{B}\big(\cup_{q=1}^{P_{b}+1}\cup_{i=1}^{p}E_{pb-p+i,q}^{\varepsilon}\big)\cup H;\nonumber
  \end{align}
  \item $T^{*}(U)=\inf\{t\in[0,T(U))|\eta^{t}(U)\in F_{\frac{\varepsilon}{2}}\cap\partial\mathcal{A}\}$ for $U\in\partial\mathcal{A}\cap (H^2)^N$.
\end{itemize}

The notations are symmetric versions of the ones in Section 3, and the difference is that we restrict the sign-changing condition for the sake of componentwise permutation.

As we proved in Section 3, the continuity of the arriving time $T^* (U)$ holds. Besides, the invariance of $T^* (U)$ is easy to check.
\begin{lemma}
$T^{*}(U)$ is continuous and invariant under the permutation $\sigma$.
\end{lemma}

To compute the $\mathbb Z_p$-index, we will use an idea from \cite{LW1}. Nevertheless, it should be noticed that the simplex in \cite{LW1} is different from ours. In order to make it work, we need to enlarge the previous set $\psi\Big(\mathbb{C}^{K\big(\sum_{b=1}^{B}P_{b}+B\big)}\Big)$. Let us select a $\sigma$-invariant set $G$ which contains $\psi\Big(\mathbb{C}^{K\big(\sum_{b=1}^{B}P_{b}+B\big)}\Big)$ by the following
\begin{align}
G=&\Bigg\{\Bigg(\sum_{k=1}^{K}\sum_{q=1}^{P_{1}+1}(-1)^{q+1}\alpha^{1,1}_{q,k}w_{1,q,k}(s^{1}_{q,k},r),\dots,\sum_{k=1}^{K}\sum_{q=1}^{P_{1}+1} (-1)^{q+1}\alpha_{q,k}^{1,p}w_{1,q,k}\Big(\frac{2\pi (p-1)}{p}+s^{1}_{q,k},r\Big),\nonumber\\
&\qquad\qquad\dots\dots\qquad\qquad\dots\dots\qquad\qquad\dots\dots\qquad\qquad\dots\dots, \nonumber\\
&\sum_{k=1}^{K}\sum_{q=1}^{P_{B}+1}(-1)^{q+1}\alpha^{B,1}_{q,k}w_{B,q,k}(s^{B}_{q,k},r),\dots,\sum_{k=1}^{K}\sum_{q=1}^{P_{B}+1} (-1)^{q+1}\alpha_{q,k}^{B,p}w_{B,q,k}\Big(\frac{2\pi (p-1)}{p}+s^{B}_{q,k},r\Big)\Bigg)\Bigg|\nonumber\\
&\alpha_{q,k}^{b,j}\geq0,\,\,s^{b}_{q,k}\in[0,2\pi),\,\,\mbox{for}\,\,\mbox{any}\,\,b=1,\dots,B,\,\,j=1\dots,p,\,\,q=1,\dots,P_{b}+1,\nonumber\\
&k=1,\dots,K\Bigg\}.\nonumber
\end{align}
Note that for any $t\geq0$ and $U\in G$, we have $tU\in G$. The difference between the set $G$ and the set $\psi\Big(\mathbb{C}^{K\big(\sum_{b=1}^{B}P_{b}+B\big)}\Big)$ is that in $G$ each coefficient of components are independent. Notice that the set $G$ contains no nontrivial fixed points of $\sigma$ due to the definitions of the functions $w_{b,q,k}$'s. We observe that due to the definition of $\mathcal{A}$ and  the property of the heat flow $\eta^{t}$, every half-line in $G$ starting at the origin intersects $\partial\mathcal{A}$. Moreover, the set $G\cap\partial\mathcal{A}$ is compact and $\sigma$-invariant. In particular, we denote
$$G_{0}=\big\{U\in G\big|n(U_{bp-p+i})=P_{b}\,\,for\,\,i=1,\dots,p\,\,and\,\,b=1,\dots,B\big\}.$$
This is to say that we define $G_{0}$ as the portion of $G$ whose elements do not degenerate in the sense of no drop-off of the sign-changing number, i.e., $n(U_{bp-p+i})=P_{b}$ for any $i=1,\dots,p$ and $b=1,\dots,B$. It is easy to see $G=\overline{G_{0}}$.

\subsection{Avoiding the Fixed Points}

In Section 3, we already proved that $\partial\mathcal{A}\cap D\cap(H^2)^N\backslash A_{\varepsilon}\neq\emptyset$. In this subsection we show that $\partial\mathcal{A}\cap D\backslash A_{\varepsilon}$ contains no fixed points of the permutation $\sigma$ action. The following lemma ensures that the flow line does not go through the fixed points of the permutation $\sigma$.

\begin{lemma}
For any $U\in\partial\mathcal{A}\cap D\backslash A_{\varepsilon}$, the flow line $\{\eta^{t}(U)\}_{t\geq0}$ contains no fixed point of the permutation $\sigma$.
\end{lemma}
\noindent{\bf Proof.}
We argue by contradiction. Suppose that there is a $t_{0}>0$ such that $\eta^{t_{0}}(U)$ is a fixed point of the permutation $\sigma$. Then we have $\eta^{t_{0}}(U)_{1}=\dots=\eta^{t_{0}}(U)_{p}$. Due to the uniqueness of the solution, $\eta^{t}(U)_{1}=\dots=\eta^{t}(U)_{p}$ for any $t\geq t_{0}$. Multiplying $u_{1}^{3}$ on both sides of the first equation of Problem (\ref{e:A14}) and integrating over $\W$, we get
\begin{align}
\frac{d}{dt}\int u_{1}^{4}+C\|u_{1}^{2}\|^{2}&= C\Bigg(\mu_{1}\int u_{1}^{6}+\sum_{i=2}^{N}\beta_{i1}\int u_{i}^{2}u_{1}^{4}\Bigg)\nonumber\\
&=C\Bigg(\big(\mu_{1}+\sum_{i=2}^{p}\beta_{i1}\big)\int u_{1}^{6}+\beta_{i1}\sum_{i=p+1}^{N}\int u_{i}^{2}u_{1}^{4}\Bigg).\nonumber
\end{align}
Combining with $\mu_{1}+\sum_{i=2}^{p}\beta_{i1}\leq0$ (assumption $(D)$ of Theorem \ref{t:main}) and the Sobolev's embedding, we have
$$\frac{d}{dt}\int u_{1}^{4}\leq-C\|u_{1}^{2}\|^{2}\leq-C|u_{1}|_{4}^{4}.$$
Hence, $\int u_{1}^4\leq Ce^{-Ct}$ follows. Therefore, for some $T_{0}>0$ and $q=1,\dots,P_{1}$, $\eta^{t}(U)\in E_{1,q}^{\frac{\varepsilon}{2}}$ for $t>T_{0}$. This is a contradiction with the fact that $U\in\partial\mathcal{A}\cap D\backslash A_{\varepsilon}$.
\begin{flushright}
\qedsymbol
\end{flushright}

\br
In fact, we can do the same computation to the other components, therefore we will have $\eta^{t}(U)\to\theta$ in $(L^{4})^{N}$.
\er

\begin{corollary}
The set $\partial\mathcal{A}\cap D\backslash A_{\varepsilon}$ contains no fixed point.
\end{corollary}

\subsection{Construction of $\sigma$-Symmetric Sets of Functions of Prescribed Nodal Numbers with Arbitrarily Large Genus}

The aim of this subsection is to prove that for any integer $k>0$, there is a compact subset $B_{k}\subset\partial\mathcal{A}\cap D\backslash A_{\varepsilon}$ satisfies $\sigma(B_{k})=B_{k}$ and $\gamma(B_{k})\geq k$. To do this, we only need to check that for the set $G$ constructed in the last subsection 4.2, it holds $\gamma(\partial\mathcal{A}\cap G \backslash A_{\varepsilon})\geq K$ since $K$ can be chosen to be arbitrarily large.

\begin{lemma}\label{c:genusG}
$\gamma(G\cap\partial\mathcal{A})=K\big(\sum_{b=1}^{B}P_{b}+B\big)$.
\end{lemma}
\noindent{\bf Proof.}
It is obvious that $\psi\Big(\mathbb{C}^{K\big(\sum_{b=1}^{B}P_{b}+B\big)}\Big)\subset G$. Hence, we have
\begin{align}
\gamma(G\cap\partial\mathcal{A})\geq\gamma\Bigg(\psi\Big(\mathbb{C}^{K\big(\sum_{b=1}^{B}P_{b}+B\big)}\Big)\cap \partial\mathcal{A}\Bigg)= K\left(\sum_{b=1}^{B}P_{b}+B\right),\nonumber
\end{align}
where the equality holds due to Borsuk's theorem. To obtain the reversed inequality, we note that the set $G$ is homeomorphic to the following subset $X$ of $\mathbb{C}^{pK\big(\sum_{b=1}^{B}P_{b}+B\big)}$:
\begin{align}
X= & \Bigg\{(z^{b,j}_{q,k})\in\mathbb{C}^{pK\big(\sum_{b=1}^{B}P_{b}+B\big)} \Bigg| arc(z^{b,1}_{q,k})=arc(z^{b,2}_{q,k})= \dots=arc(z^{b,p}_{q,k}) \nonumber\\
& \mbox{for}\,\,\mbox{any}\,\,b=1,\dots,B,\,\,q=1,\dots,P_{b},\,\,k=1,\dots,K\Bigg\}.\nonumber
\end{align}
In fact we may define $\xi : G\to X$ by
\begin{align}
\xi:\Bigg(\sum_{k=1}^{K}\sum_{q=1}^{P_{b}+1} (-1)^{q+1}\alpha_{q,k}^{b,j}w_{b,q,k}\Big(\frac{2\pi (j-1)}{p}+s^{b}_{q,k},\cdot\Big)\Bigg) \to \Bigg(e^{is^{b}_{q,k}}\alpha_{q,k}^{b,j} \Bigg).\nonumber
\end{align}
 To distinguish between the spaces $\mathbb{C}^{K\big(\sum_{b=1}^{B}P_{b}+B\big)}$ and $\mathbb{C}^{pK\big(\sum_{b=1}^{B}P_{b}+B\big)}$, we denote their vectors by $(z_{b,q,k})$ and $(z_{q,k}^{b,j})$ respectively. On the other hand, we have a continuous map $f:X\to\mathbb{C}^{K\big(\sum_{b=1}^{B}P_{b}+B\big)}$ written as
\begin{align}
f(z_{q,k}^{b,j})=\Big(\sum_{j=1}^{p}z_{q,k}^{b,j}\Big).\nonumber
\end{align}
Notice that $f^{-1}(0)=0$ and $f(e^{\frac{2\pi i}{p}}z_{q,k}^{b,j})=e^{\frac{2\pi i}{p}}f(z_{q,k}^{b,j})$. The reversed inequality follows from the identity $f\circ\xi(\sigma U)=e^{\frac{2\pi i}{p}}f\circ\xi(U)$.
\begin{flushright}
\qedsymbol
\end{flushright}

By Proposition \ref{p:L4} for $\varepsilon>0$ small enough, $E_{bp-p+i,q}^{\varepsilon}$ is invariant under the heat flow $\eta^{t}(\cdot)$ for any  $b=1,\dots,B$, $i=1,\dots,p$ and $q=1,\dots,P_{b}$. For the sake of convenience, we will denote the set $E_{bp-p+i,q}^{\varepsilon}$ by $E_{j,q}^{\varepsilon}$ with $b=1,\dots,B$, $i=1,\dots,p$, $q=1,\dots,P_{b}$ and $j=pb-p+i$. The complete invariant set of $E_{j,q}^{\varepsilon}$'s defined by
$$C(E_{j,q}^{\varepsilon})=\{U=(u_{1},\dots,u_{N})\in(H_{r}^{2})^{N}|\exists t_{0}\geq0\,\,s.t.\,\,\eta^{t_{0}}(U)\in E_{j,q}^{\varepsilon}\},$$
for any admissible $j$'s and $q$'s. The concept of complete invariant set is from \cite{LS}. As we did in Section 3.2, define a continuous cut-off function $\phi:[0,\infty)\to[0,1]$ by
$\phi(s)=1$ for $s\geq \varepsilon$, $\phi(s)=0$ for $s\leq\frac{\varepsilon}{2}$, and $\phi(s)=\frac{2s}{\varepsilon}-1 $ for  $s\in (\frac{\varepsilon}{2},\varepsilon)$.
We Denote $$A_{\varepsilon}=\cup_{b=1}^{B}\cup_{i=1}^{p}\cup_{q=1}^{P_{b}+1}C(E_{bp-p+i,q}^{\varepsilon})\cup H.$$

 We will prove $\gamma(\partial\mathcal{A}\cap G\backslash A_{\varepsilon})\geq K - (p-1)(\sum_{b=1}^{B}P_{b}+B)$, which completes the proof. To obtain this, we define a mapping $h: G\cap A_{\varepsilon}\to G$ as
$$h(U)=\Big(\sum_{q=1}^{P_{1}+1}\phi\big(|\eta^{T^{*}(U)}(U)_{1,q}|_{4})U_{1,q},\dots, \sum_{q=1}^{P_{B}+1}\phi\big(|\eta^{T^{*}(U)}(U)_{N,q}|_{4}\big)U_{N,q}\Big)$$
for $U\in G_{0}$ and for $U\in G\backslash G_{0}$, let
$$h(U)=\Big(\sum_{q=1}^{P_{1}+1}\phi\big(|U_{1}\cdot\chi_{\W_{1,q}}|_{4})U_{1}\cdot\chi_{\W_{1,q}},\dots, \sum_{q=1}^{P_{B}+1}\phi\big(|U_{N}\cdot\chi_{\W_{N,q}}|_{4}\big)U_{N}\cdot\chi_{\W_{N,q}}\Big).$$
Here, $\phi$ is the cut-off function in \eqref{cutoff}, $\chi_{A}$ is the characteristic function of the set $A\subset\W$. In fact, for $U\in G$, $U_{bp-p+i,q}=U_{bp-p+i}\cdot\chi_{\W_{b,q}}$ for $q=1\dots,P_{b}+1$, $i=1,\dots,p$ and $b=1\dots,B$. And $n(U_{bp-p+i})\leq P_{b}$ for $i=1,\dots,p$ and $b=1,\dots,B$. The "$\leq$"'s hold strictly for at least one of admissible $(b,j)$'s. It is easy to see that the map $h$ is continuous. Then, we have the following claim.
\begin{lemma}\label{c:phi}
For any $U\in\partial\mathcal{A}\cap G_{0}\cap A_{\varepsilon}$, there are admissible couples $(j_{1},q_{1})$ and $(j_{2},q_{2})$ such that $\phi(|\eta^{T^{*}(U)}(U)_{j_{1},q_{1}}|_{4})=0$ and $\phi(|\eta^{T^{*}(U)}(U)_{j_{2},q_{2}}|_{4})>0$.
\end{lemma}
\noindent{\bf Proof.}
We notice that $\eta^{t}(U)$ always stays on $\partial\mathcal{A}$, which implies that $\|\eta^{t}(U)\|\geq C>0$ for any $t>0$. Let us assume that $\phi(|\eta^{T^{*}(U)}(U)_{j,q}|_{4})=0$ for any $(j,q)$ admissible. This gives $\sum_{j=1}^{N}\int|u_{j}|^{4}<C\varepsilon^4$ as $t$ goes large. Multiplying $u_{j}$ on the both sides of the $j$-th equation of Problem (\ref{e:A14}) and summing up with respect to $j$, we have
\begin{align}
\frac{1}{2}\frac{\partial}{\partial t}|U(t)|^{2}_{2}+\|U(t)\|^{2}\leq\sum_{j=1}^{N}\Bigg(\mu_{j}\int|u_{j}(t)|^4 +\sum_{i\neq j}\beta_{ij}\int u_{i}(t)^{2}u_{j}(t)^{2}\Bigg)<C\varepsilon\nonumber
\end{align}
when $t$ is large. Using the openness of $\mathcal{A}$ in $(H^1)^N$ and the invariance of $\partial\mathcal{A}$, we have $\frac{\partial}{\partial t}|U(t)|_{2}^{2}\leq-C$ for $t>0$. This is a contradiction, and the proof is complete.
\begin{flushright}
\qedsymbol
\end{flushright}

\br
The lemma implies that for any $U\in\partial\mathcal{A}\cap G_{0}\cap A_{\varepsilon}$, there are two admissible couples $(j_{1},q_{1})$ and $(j_{2},q_{2})$ such that $h(U)_{j_{1},q_{1}}=0$ and $h(U)_{j_{2},q_{2}}\neq0$.
\er

\begin{lemma}\label{ge}
 It holds $\gamma(\partial\mathcal{A}\cap D\backslash A_{\varepsilon})=\infty$.
\end{lemma}
\noindent{\bf Proof.}
We use the notations in the proof of Lemma \ref{c:genusG}. To proceed our computation, we need an upper bound of $\gamma\left(h(\partial\mathcal{A}\cap G\cap A_{\varepsilon})\right)$. Due to above deduction, for any $U\in h(\partial\mathcal{A}\cap G\cap A_{\varepsilon})$, we can find a admissible couple $(j,q)$ such that $|U_{j,q}|_{4}=0$. Translating into the version in $G$, we have that there are some $b_{0}=1,\dots,B$, $j_{0}=1,\dots,p$ and $q_{0}=1,\dots,P_{b_{0}}$ such that $|z_{q_{0},k}^{b_{0},j_{0}}|=0$ for any $k=1,\dots,K$. Instead of estimating $\gamma\left(h(\partial\mathcal{A}\cap G_{0}\cap A_{\varepsilon})\right)$, we will give an upper bound of $\gamma\left(\xi\circ h(\partial\mathcal{A}\cap G_{0}\cap A_{\varepsilon})\right)$ since $\gamma\left(h(\partial\mathcal{A}\cap G\cap A_{\varepsilon})\right)\leq \gamma\left(\xi\circ h(\partial\mathcal{A}\cap G_{0}\cap A_{\varepsilon})\right)$.

We divide $h(\partial\mathcal{A}\cap G_{0}\cap A_{\varepsilon})$ into two parts. The first part is defined by
\begin{align}
L=&\Bigg\{(z_{q,k}^{b,j})\in\xi\circ h(\partial\mathcal{A}\cap G\cap A_{\varepsilon})\Bigg|\exists b_{0}=1,\dots,B,~~\exists q_{0}=1,\dots, P_{b_{0}}\nonumber\\
&s.t.~~|z_{q_{0},k}^{b_{0},j}|=0~~\forall j=1,\dots,p,~~\forall k=1,\dots,K\Bigg\}.\nonumber
\end{align}

It is easy to see that for any $(z_{b,q,k})\in f(L)\subset\mathbb{C}^{K(\sum_{b=1}^{B}P_{b}+B)}$, we can find some $b_{0}=1,\dots,B$ and $q_{0}=1,\dots, P_{b_{0}}$ such that for any $k=1,\dots,K$, $z_{b_{0},q_{0},k}=0$. We define a subspace $W\subset\mathbb{C}^{K\big(\sum_{b=1}^{B}P_{b}+B\big)}$ by
\begin{align}
W=&\Big\{(z_{b,q,k})\in\mathbb{C}^{K\big(\sum_{b=1}^{B}P_{b}+B\big)}|z_{1,1,k}=\dots=z_{1,P_{1}+1,k}=\dots=z_{B,1,k}=\dots\nonumber\\
&\dots= z_{B,P_{B}+1,k}\,\,\forall k=1,\dots,K\Big\}.\nonumber
\end{align}
Obviously, $\psi(W)$ is invariant under the permutation $\sigma$. Since for every element $U\in\psi(W)$, $n(U_{bp-p+i})=P_{b}$ for $i=1,\dots,p$ and $b=1,\dots,B$, we have $\psi(W)\cap \partial\mathcal{A}\subset G_{0}\cap\partial\mathcal{A}$. To continue the proof, we need the following lemma:

\begin{lemma}\label{l:homotopynbh}
There is an $\varepsilon_{0}>0$ such that for any $a=(a_{1},\dots,a_{K(\sum_{b=1}^{B}P_{b}+B)})\in f(L)$, $\sum_{l=1}^{K(\sum_{b=1}^{B}P_{b}+B)}|a_{l}|\geq\varepsilon_{0}$.
\end{lemma}
\noindent{\bf Proof.}
We argue it by contradiction. Then there is a sequence $a^{(n)} \in f(L)$ such that $|a^{(n)}|\to 0$, which implies
that there is a sequence $U^{(n)}\in \partial\mathcal{A}\cap G\cap A_{\varepsilon}$ such that $\sum_{j=1}^{N}|u_{j}^{(n)}|_{4}\to 0$. Using the fact that $T(U)=\infty$ for $U\in\partial\mathcal{A}$ and similar computations in Proposition \ref{p:L4} and Lemma \ref{c:phi}, we will have a contradiction.
\begin{flushright}
\qedsymbol
\end{flushright}

We now return to the proof of Lemma \ref{ge}. On one hand, we have $dim(W)=K$.
On the other hand, due to the definition of the space $W$ and Lemma \ref{l:homotopynbh}, we have
\begin{align}
f(L)\subset\mathbb{C}^{K\big(\sum_{b=1}^{B}P_{b}+B\big)}\backslash W_{\delta}\overset{\sigma}{\simeq}\mathbb{S}^{K\big(\sum_{b=1}^{B}P_{b}+B-1\big)-1}.\nonumber
\end{align}
Here, $\delta>0$ is small, $W_{\delta}$ represents the $\delta$-neighbourhood of $W$, and the symbol $\overset{\sigma}{\simeq}$ means that two topology spaces are homotopy equivalent via a homotopy $F$ satisfies $F(t,e^{\frac{2\pi i}{p}}z)=e^{\frac{2\pi i}{p}}F(t,z)$ and $\mathbb{S}^{m-1}$ denote the unit sphere in $\mathbb{C}^{m}$. Hence, we have
\begin{align}
\gamma(L)\leq\gamma(f(L))\leq K\Bigg(\sum_{b=1}^{B}P_{b}+B-1\Bigg).\nonumber
\end{align}

Now we give an estimate on another part of $\xi \circ h(\partial\mathcal{A}\cap G_{0}\cap A_{\varepsilon})$. Define the set
$$M:=\xi\circ h(\partial\mathcal{A}\cap G\cap A_{\varepsilon})\backslash L_{\delta},$$
where $\delta>0$ is small. According to the previous deduction, the elements $(z_{q,k}^{b,j})$ in $M$ satisfy that there are $b_{0}=1,\dots,B$ and $q_{0}=1,\dots,P_{b_{0}}+1$ such that there are two $j_{1},~j_{2}=1,\dots,p$ with $j_1\neq j_2$ and
\begin{itemize}
  \item $|z_{q_{0},k}^{b_{0},j_{1}}|=0$ for any $k=1,\dots,K$;
  \item $\sum_{k=1}^{K}|z_{q_{0},k}^{b_{0},j_{2}}|>0$.
\end{itemize}
Using these for any $b=1\dots,B$ and $q=1,\dots,P_{b_{1}}+1$ we define
\begin{align}
N_{q}^{b}=\Bigg\{(z_{q,k}^{b,j})\in M\Bigg| \mbox{there exist $j_1\neq j_2$ such that}\; \sum_{k=1}^{K} \left| z_{q,k}^{b,j_1}\right|=0, \sum_{k=1}^{K}|z_{q,k}^{b,j_{2}}|>0\Bigg\}.\nonumber
\end{align} The following two properties follow immediately
\begin{itemize}
  \item $\cap_{i=0}^{p-1}\sigma^{i}\big(\xi^{-1}(N_{q}^{b})\big)=\emptyset$ for any $b$ and $q$;
  \item $\xi^{-1}(M)\subset\cup_{b=1}^{B}\cup_{q=1}^{P_{b}+1}\cup_{i=0}^{p-1}\sigma^{i}\big(\xi^{-1}(N_{q}^{b})\big)$.
\end{itemize}
Therefore, using the last property of the index $\gamma$ in Section 4.1 we have
\begin{align}
\gamma(M)=\gamma\big(\xi^{-1}(M)\big)\leq\sum_{b=1}^{B}\sum_{q=1}^{P_{b}+1}\gamma\left(\mathbb{Z}_{p}(\xi^{-1}(N_{q}^{b}))\right) \leq(p-1)\left(\sum_{b=1}^{B}P_{b}+B\right).\nonumber
\end{align}
Now we can compute
\begin{align}
\gamma\left(\overline{\partial\mathcal{A}\cap G\cap A_{\varepsilon}}\right)
&\leq\gamma\left(\overline{h(\partial\mathcal{A}\cap G\cap A_{\varepsilon})}\right)\leq \gamma\left(\overline{\xi\circ h(\partial\mathcal{A}\cap G\cap A_{\varepsilon})}\right)\nonumber\\
&\leq\gamma(L_{\delta})+\gamma(M_{\delta})=\gamma(L)+\gamma(M)\nonumber\\
&\leq K\Bigg(\sum_{b=1}^{B}P_{b}+B-1\Bigg)+(p-1)\left(\sum_{b=1}^{B}P_{b}+B\right).\nonumber
\end{align}
The above result implies that
\begin{align}
&\;\;\;\;\;\gamma\big((\partial\mathcal{A}\backslash A_{\varepsilon})\cap G\big)\\
&\geq\gamma(G\cap\partial\mathcal{A})- \gamma\big(\overline{A_{\varepsilon}\cap\partial\mathcal{A}\cap G}\big)\nonumber\\
&\geq K\Bigg(\sum_{b=1}^{B}P_{b}+B\Bigg)-K\Bigg(\sum_{b=1}^{B}P_{b}+B-1\Bigg)-(p-1)\left(\sum_{b=1}^{B}P_{b}+B\right)\nonumber\\
&=K-(p-1)\left(\sum_{b=1}^{B}P_{b}+B\right).\nonumber
\end{align}
The proof of Lemma \ref{ge} is complete since $K$ is arbitrarily large.
\begin{flushright}
\qedsymbol
\end{flushright}

\br
We note that
$\partial\mathcal{A}\cap D\backslash A_{\varepsilon}$ is an invariant set of the heat flow, from which, a sequence of compact sets with unbounded genus can be selected.
\er

\subsection{The Existence of Multiple Equilibria Having the Same Componentwise Prescribed Nodes}
In this subsection, we complete the proof of the main result Theorem \ref{t:main}.

\begin{lemma}\label{def}
 Let $c\in J(\partial\mathcal{A}\cap D\backslash A_{\varepsilon})$. If there are positive numbers $\alpha$ and $\varepsilon$ such that for any $U\in J^{-1}[c-\varepsilon,c+\varepsilon]\cap\partial\mathcal{A}\cap D\backslash A_{\varepsilon}$, $|\Delta u_{j}-\lambda u_{j}+\mu u_{j}^{3}+\sum_{i\neq j}\beta u_{i}^{2}u_{j}|_{2}\geq\alpha$ for some $j=1,\dots,N$, then there is $T>0$ independent of $U$ such that $\eta^{T}(U)\in J^{c-\varepsilon}\cap\partial\mathcal{A}\cap D\backslash A_{\varepsilon}$.
\end{lemma}
\noindent{\bf Proof.}
Let $T=\frac{4\varepsilon}{\alpha^2}>0$. Notice that $\partial\mathcal{A}\cap D\backslash A_{\varepsilon}$ is invariant under the flow $\eta$. If $\eta^{T}(U)\in J^{c-\varepsilon}\cap\partial\mathcal{A}\cap D\backslash A_{\varepsilon}$, the proof is complete. Otherwise, assume $J(\eta^{T}(U))> c-\varepsilon$. Then $\eta^{t}(U)\in J^{-1}[c-\varepsilon,c+\varepsilon]\cap\partial\mathcal{A}\cap D\backslash A_{\varepsilon}$ for any $t\in[0,T]$. Compute that
\begin{align}
\frac{d}{dt}J(\eta^{t}(U))&=-\sum_{j=1}^{N}\int|\partial_{t}u_{j}(t,\cdot)|_{2}^{2}\nonumber\\
&=-\sum_{j=1}^{N}\int\big|\Delta u_{j}-\lambda_{j} u_{j}+\mu_{j} u_{j}^{3}+\sum_{i\neq j}\beta_{ij} u_{i}^{2}u_{j}\big|_{2}^{2}(t)\nonumber\\
&\leq-\alpha^{2}.\nonumber
\end{align}
Therefore,
\begin{align}
c-\varepsilon\leq c+\int_{0}^{T}\frac{d}{dt}J(\eta^{t}(U))dt\leq c-\alpha^{2}T=c-4\varepsilon.\nonumber
\end{align}
This is a contradiction.
\begin{flushright}
\qedsymbol
\end{flushright}

Define
$$\Gamma_{k}=\{A\subset\partial\mathcal{A}\cap D\backslash A_{\varepsilon}| \mbox{$A$ is $\sigma$ invariant compact set},\,\gamma(A)\geq k\}.$$
By Lemma \ref{ge}, $\Gamma_{k}\neq\emptyset$ for large $k$, and the values
$$c_{k}=\inf_{A\in\Gamma_{k}}\sup_{u\in A}J(U)$$
are well-defined. Using Lemma \ref{def} and some classical arguments as \cite[Propersition 8.5]{R} it is easy to verify the following.
\begin{lemma}
(i). $K_{c_{k}}\cap\partial\mathcal{A}\cap D\backslash A_{\varepsilon}\neq\emptyset$ for large $k$.

(ii). If $c:=c_{j}=\dots=c_{j+l}$, then $\gamma(K_{c}\cap\partial\mathcal{A}\cap D\backslash A_{\varepsilon})\geq l+1$.
\end{lemma}
\noindent{\bf Proof.}
A standard argument ensures that there is a sequence $U_{n}\in\partial\mathcal{A}\cap D\backslash A_{\varepsilon}$, such that
\begin{itemize}
  \item $J(U_{n})\to c_{k}$;
  \item for any $j=1,\dots,N$, $\Delta u_{n,j}-\lambda_{j} u_{n,j}+\mu_{j} u_{n,j}^{3}+\sum_{i\neq j}\beta_{ij} u_{n,i}^{2}u_{n,j}\to0$ in $L^{2}$,
\end{itemize}
where $u_{n,j}$ is the $j$-th component of $U_{n}$. The second assertion implies that $\nabla_{u_{j}}J(U_{n})\to\theta$ in $H^{-1}$. Then $U_{n}\to U$ in $(H^1)^N$ for some $U$ since the energy functional $J$ satisfies the (PS) condition. It is easy to see that $U$ is a critical point of the energy $J$. Therefore, $U$ is of class $(H^2)^N$ due to the elliptic regularity theory. Now we prove that $U\in\partial\mathcal{A}\cap D\backslash A_{\varepsilon}$.

We have $U\in\partial\mathcal{A}\setminus A_{\varepsilon}$ due to the elliptic regularity and $(L^4)^N$-norm is continuous in $(H^1)^N$.

To show $U\in D$, denote $U= (u_1, ..., u_N)$. $U_{n}$'s are continuous since they are of class $(H^2)^N$.
We have for $b=1,...,B$, $n\big((U_{n})_{pb-p+i}\big)=P_{b}$ for $i=1,...,p$, and $\big|(U_{n})_{pb-p+i,q}\big|_{4}\geq\frac{\varepsilon}{2}$ for
$i=1,...,p$ and $q=0,\dots,P_{b}$. Here $(U_{n})_{pb-p+i,q}$ denotes the $q$-th bump of the $(pb-p+i)$-th component of $U_{n}$.

Since $U_{n}\to U$ in $(H^1)^N$, for $n$ large, we have $\big|(U_{n})_{pb-p+i}-u_{pb-p+i}\big|_{4}<\frac{\varepsilon}{100}$.
 Therefore, for fixed $b=1,...,B, i=1,...,p$,
we can find a sequence $x_{0},\dots,x_{P_{b}}\in\Omega$ such that
\begin{itemize}
  \item $0<|x_{0}|<\dots<|x_{P_{b}}|<\infty$;
  \item $u_{pb-p+i}(x_{k})\cdot u_{pb-p+i}(x_{k+1})<0$ for $k=0,\dots,P_{b}-1$.
\end{itemize}
Therefore, $U\in D$. Assertion (i) is proved.

Assertion (ii) can be proved by some arguments for the classical genus \cite[Propersition 8.5]{R}, and we omit it here.
\begin{flushright}
\qedsymbol
\end{flushright}

\section*{Acknowledgement}
Li would like to express his sincere appreciation to Thomas Bartsch, Marek Fila, Desheng Li and Pavol Quittner
for their useful communication and comments. Li especially thanks Pavol Quittner for his
kindly offering \cite{
Q1,Q2,Q2003,Q3}. Wang thanks Zhaoli Liu for some useful discussion. This work is supported by NSFC (11771324, 11831009).

\Vs\Vs


{\footnotesize

\begin {thebibliography}{44}

\bibitem{AB}Ackermann, N., Bartsch, T., Superstable Manifolds of Semilinear Parabolic Problems. Journal of Dynamics and Differential Equations, 17 (2005) 115-173.
\bibitem{AA}Akhmediev, N., Ankiewicz, A., Partially Coherent Solitons on a Finite Background. Physical Review Letters, 82 (1999) 2661-2664.
\bibitem{Amann87}Amann, H., On abstract parabolic fundamental solutions. Journal of the Mathematical Society of Japan, 39 (1987) 93-116.
\bibitem{Amann}Amann, H., Linear and Quasilinear Parabolic Problems, Birkh\"{a}user Verlag, Basel, 1995.
\bibitem{AC}Ambrosetti, A., Colorado, E., Standing waves of some coupled nonlinear Schr\"{o}dinger equations. Journal of the London Mathematical Society, 75 (2007) 67-82.
\bibitem{AR}Ambrosetti, A., Rabinowitz, P., H., Dual variational methods in critical point theory and applications. Journal of Functional Analysis, 14 (1973) 349-381.

\bibitem{An}Angenent, S., Nodal properties of solutions of parabolic equations. Rocky Mountain J.
Math., 21 (1991) 585-592.

\bibitem{AWY} Ao, W., Wei, J., Yao, W., Uniqueness and non-degeneracy of sign-changing radial solutions to an almost critical elliptic problem, arXiv:1510.04678v1.

\bibitem{BW}Bartsch, T., Wang, Z.-Q., Note on ground states of nonlinear Schr\"{o}dinger systems. Partial Differential Equations. 3 (2006) 200-207.


\bibitem{BDW} Bartsch, T., Dancer, E.N., Wang, Z.-Q.,
     A Liouville theorem, a-priori bounds, and bifurcating branches of positive solutions for a nonlinear elliptic system.
    {Calc. Var. Partial Differential Equations} {37} (2010) 345-361.

\bibitem{BWillem} Bartsch, T., Willem, M., Infinitely many radial solutions of a semilinear elliptic problem on $\mathbb R^N$. Archive for Rational Mechanics and Analysis, 124 (1993) 261-276.
\bibitem{BWW}Bartsch, T., Wang, Z.-Q., Willem, M., The Dirichlet problem for superlinear elliptic equations. Handbook of Differential Equations: Stationary Partial Differential Equations. Vol 2, 1-55, Elsevier Science and Technology, 2005.


\bibitem{CL}Cazenave, T., Lions, P, L., Solutions globales d'��quations de la chaleur semi lin��aires. Communications in Partial Differential Equations, 9 (1984), 955-978.

\bibitem{Ch0}Chang, K.-C., Infinite Dimensional Morse Theory and Multiple Solutions Problems.
Birkh\"{a}user Boston, 1993.


\bibitem{Ch}Chang, K.-C., Heat method in nonlinear elliptic equations. Topological methods, variational methods and their applications (Tianyuan, 2002), 65-76, World Sci., River Edge, NJ, 2003.
\bibitem{CMT}Conti, M., Merizzi, L., Terracini, S., Radial Solutions of Superlinear Equations on $\mathbb{R}^{N}$. Part I: Global Variational Approach. Archive for Rational Mechanics and Analysis, 153 (2000) 291-316.
\bibitem{DWW}Dancer, E, N., Wei, J., Weth, T., A priori bounds versus multiple existence of positive solutions for a nonlinear Schr\"{o}dinger system. Annales De Linstitut Henri Poincare, 27 (2010) 953-969.

\bibitem{DM}Daners, D., Medina, P, K., Abstract Evolution Equations, Periodic Problems, and Applications. Longman Scientific and Technical, 1992.


\bibitem{FL}Fila, M., Levine, H. A., On the Boundedness of Global Solutions of Abstract Semilinear Parabolic Equations. Journal of Mathematical Analysis and Applications, 216 (1997) 654-666.

\bibitem{G}Galaktionov, V., Geometric sturmian theory of nonlinear parabolic equations and applications. Chapman and Hall/CRC Florida, 2004.
\bibitem{GNN}Gidas, B., Ni, W, M., Nirenberg, L., Symmetry and related properties via the maximum principle. Communications in Mathematical Physics, 68 (1979) 209-243.
\bibitem{H}Henry, D., Geometric Theory of Semilinear Parabolic Equations. Lecture Notes in Mathematics, 840, Springer-Verlag, Berlin, 2008.
\bibitem{K}Kwong, M, K., Uniqueness of positive solutions of $\Delta u-u+u^{p}=0$ in $\mathbb{R}^n$. Archive for Rational Mechanics and Analysis, 105 (1989) 243-266.
\bibitem{LM1111111}Li, H., Meng, L., Notes on A Superlinear Elliptic Problem, preprint.
\bibitem{LWei1}{Lin, T.-C., Wei, J.,}
{Ground state of $N$ coupled nonlinear Schr\"odinger equations
in $ R^n, n\leq3$.}
 {Comm. Math. Phys.} {255} (2005) 629--653.

\bibitem{LWei2}
    {Lin, T.-C., Wei, J.,}
    { Spikes in two coupled nonlinear Schr\"{o}dinger equations.}
    { Ann. Inst. H. Poincar\'e Anal. Non Lin\'eaire } {22} (2005) 403--439.



\bibitem{LM}Lions, J. L., Magenes, E., Non-Homogeneous Boundary Value Problems and Applications, Volum 1. Springer-Verlag, Berlin, 1972.


\bibitem{LLW}Liu, J., Liu, X., Wang, Z.-Q., Multiple mixed states of nodal solutions for nonlinear Schr\"{o}dinger systems. Calculus of Variations and Partial Differential Equations, 52 (2015) 565-586.

\bibitem{LS}Liu, Z., Sun, J., Invariant Sets of Descending Flow in Critical Point Theory with Applications to Nonlinear Differential Equations. Journal of Differential Equations, 172 (2001) 257-299.
\bibitem{LW1}Liu, Z., Wang, Z.-Q.. Multiple bound states of nonlinear Schr\"{o}dinger systems. Comm. Math. Phys., 282 (2008) 721-731.
\bibitem{LW}Liu, Z., Wang, Z.-Q., Ground States and Bound States of a Nonlinear Schr\"{o}dinger System. Advanced Nonlinear Studies, 10 (2010) 175-193.
\bibitem{LW2}Liu, Z., Wang, Z.-Q., Vector solutions with prescribed component-wise nodes for a Schr\"{o}dinger system. Anal. Theory Appl., 35 (2019) 288-311.


\bibitem{Lu}Lunardi, A., Analytic semigroups and optimal regularity in parabolic problems. Progress in Nonlinear Differential Equations and their Applications, 16. Birkh\"{a}user Verlag, Basel, 1995.


\bibitem{Mat}Matano, H., Nonincrease of the lap number of a solution for a one-dimensional semi-linear parabolic equation, J. Fac. Sci. Univ. Tokyo Sect. 1A Math., 29 (1982) 401-441.

\bibitem{MS}Mitchell, M., Segev, M., Self-trapping of incoherent white light beams. Nature, 387 (1997) 880-883.

\bibitem{NTTV}
{Noris, B., Tavares, H., Terracini, S., Verzini, G.,}
{Uniform H\"older bounds for nonlinear Schr\"odinger systems with strong competition.}
 {Comm. Pure and Appl. Math.} {63} (2010) 267-302.



\bibitem{PrS}Pr\"uss, J., Sohr, H., Imaginary powers of elliptic second order differential operators in $L_p$-spaces. Hiroshima Mathematical Journal. 23 (1993) 395-418.




\bibitem{Q1}Quittner, P., Boundedness of trajectories of parabolic equations and stationary solutions via dynamical methods. Differential and Integral Equations., 7 (1994) 1547-1556.
\bibitem{Q2}Quittner, P., Signed solutions for a semilinear elliptic problem. Differential and Integral Equations, 11 (1998) 551-559.
\bibitem{Q2003}Quittner, P., Continuity of the blow-up time and a priori bounds for solutions in superlinear parabolic problems. Houston Journal of Mathematics, 29 (2003) 757-799.
\bibitem{Q3}Quittner, P., Multiple equilibria, periodic solutions and a priori bounds for solutions in superlinear parabolic problems. Nonlinear Differential Equations and Applications Nodea, 11 (2004) 237-258.

\bibitem{Q}Quittner, P., Souplet, P., Superlinear parabolic problems. Birkh\"{a}user, Basel, 2007.



\bibitem{R}Rabinowitz, P., Minimax Methods in Critical Point Theory with Applications to Differential Equations, CBMS Regional Conf. Ser. in Math., vol. 65. American Mathematical Society, Providence, 1986.

\bibitem{S}Sirakov, B.,
    Least energy solitary waves for a system of nonlinear Schr\"odinger equations in $ \mathbb R^n$.
    {Comm. Math. Phys.} {271} (2007) 199-221.


\bibitem{See}Seeley, R., Interpolation in $L_p$ with boundary conditions. Studia Math. 44 (1972) 47-60.


\bibitem{Struwe}Struwe, M., Superlinear elliptic boundary value problems with rotational symmetry.
 Archiv der mathematik. 39 (1982) 233-240.


\bibitem{T}Tartar, B, L., An Introduction to Sobolev Spaces and Interpolation Spaces. Springer-Verlag Berlin Heidelberg, 2007.

\bibitem{Ta}Tanaka, S. Uniqueness of sign-changing radial solutions for $\Delta u-u+|u|^{p-1}u=0$ in some ball and annulus. Journal of Mathematical Analysis and Applications, 439 (2016) 154-170.

\bibitem{TV}
     {Terracini, S., Verzini, G.,}
     {Multipulse Phase in $k$-mixtures of Bose-Einstein condensates.}
     { Arch. Rat. Mech. Anal.} {194} (2009) 717-741.




\bibitem{TW}Tian, R., Wang, Z.-Q., Multiple solitary wave solutions of nonlinear Schr\"{o}dinger systems. Topo. Methods Nonlinear Anal. 37 (2011) 203-223.
\bibitem{Wang}Wang, Z.-Q.,
A $\mathbb Z_p$-Borsuk-Ulam Theorem. {Chinese Bulletin of
Science}, {34} (1989) 1153-1157.

\bibitem{WW}Wei, J., Weth, T., Radial Solutions and Phase Separation in a System of Two Coupled Schr\"{o}dinger Equations. Archive for Rational Mechanics and Analysis, 190 (2008) 83-106.
\bibitem{W}Willem, M., Minimax Theorems. Birkh\"{a}user Boston, 1996.

\end {thebibliography}
}
\end{document}